\documentclass{IEEEtran4PSCC}
\usepackage{balance}
\usepackage{stfloats}

\ifCLASSINFOpdf
   \usepackage[pdftex]{graphicx}
\else
   \usepackage[dvips]{graphicx}
\fi
\usepackage[cmex10]{amsmath}
\usepackage{wasysym}
\usepackage{cite}
\usepackage{thmtools}
\usepackage{amsthm}
\usepackage{graphicx}
\usepackage{multirow}
\usepackage[export]{adjustbox}
\usepackage[normalem]{ulem}
\usepackage{amsmath}
\usepackage{amssymb}
\usepackage[bb=boondox]{mathalfa}
\usepackage{enumitem}
\usepackage{pgfplots}
\pgfplotsset{compat=1.9}
\usepackage{booktabs}
\usepackage{hyperref}
\usepackage{pdfpages}

\ifCLASSOPTIONcompsoc
    \usepackage[caption=false, font=normalsize, labelfont=sf, textfont=sf]{subfig}
\else
\usepackage[caption=false, font=footnotesize]{subfig}
\fi

\usepackage{tikz}
\usetikzlibrary{decorations.pathmorphing}
\usetikzlibrary{decorations.pathreplacing}
\usetikzlibrary{shapes}

\newtheoremstyle{remark*}
  {\topsep}
  {\topsep}
  {}
  {}
  {\itshape}
  {:}
  {.5em}
  {\thmname{#1}\thmnumber{ #2}\thmnote{ (#3)}}
\theoremstyle{remark*}
\newtheorem*{remark*}{Remark}

\newcommand{\fat}[1]{\boldsymbol{#1}}

\newcommand{\braceit}[1]{\left({#1}\right)}

\usepackage{mathtools}
\newcommand\norm[1]{\left\lVert#1\right\rVert}

\usepackage{esvect}

\DeclareMathSymbol{\shortminus}{\mathbin}{AMSa}{"39}
\def\shortplus{
\begin{tikzpicture}
    \node[draw=none,scale=0.5,inner sep=0] at (0,0) {+};
\end{tikzpicture}
}

\makeatletter
\newcommand{\oset}[3][0ex]{%
  \mathrel{\mathop{#3}\limits^{
    \vbox to#1{\kern-2\ex@
    \hbox{$\scriptstyle#2$}\vss}}}}
\makeatother

\def\compressortwo#1#2{
\begin{scope}[shift={#1}, rotate=#2, scale = 1.2]
    \draw[black, line width=0.25mm, fill = gray!50] (0.7,0) -- (0.7,0.2) -- (1.1,0.1) -- (1.1,-0.1) -- (0.7,-0.2) -- (0.7,0);
\end{scope}
}

\definecolor{maincolor}{HTML}{032F99}
\definecolor{blue}{HTML}{032F99}
\definecolor{red}{HTML}{e05a87}
\definecolor{green}{HTML}{228b22}

\renewcommand{\baselinestretch}{0.97}

\begin{document}
\begingroup
\allowdisplaybreaks
%
\title{Multi-Stage Linear Decision Rules for Stochastic Control of Natural Gas Networks with Linepack}


\author{\IEEEauthorblockN{
Vladimir Dvorkin\IEEEauthorrefmark{1},
Dharik Mallapragada\IEEEauthorrefmark{1}, 
Audun Botterud\IEEEauthorrefmark{1},
Jalal Kazempour\IEEEauthorrefmark{2} and 
Pierre Pinson\IEEEauthorrefmark{2}}
\IEEEauthorblockA{\IEEEauthorrefmark{1}Massachusetts Institute of Technology, Cambridge, MA, USA.}
\IEEEauthorblockA{\IEEEauthorrefmark{2}Technical University of Denmark, Kng. Lyngby, Denmark.}
\vspace{-1cm}
}

\maketitle

\begin{abstract}
The disturbances from variable and uncertain renewable generation propagate from power systems to natural gas networks, causing gas network operators to adjust gas supply nominations to ensure operational security. To alleviate expensive supply adjustments, we develop control policies to leverage instead the flexibility of linepack -- the gas stored in pipelines -- to balance stochastic gas extractions. These policies are based on multi-stage linear decision rules optimized on a finite discrete horizon to guide controllable network components, such as compressors and valves, towards feasible operations. Our approach  offers several control applications. First, it treats the linepack as a main source of flexibility to balance disturbances from power systems without substantial impacts on nominal gas supply. Second, these policies can be optimized to minimize the variability (due to intermittency of renewables) and variance (due to their uncertainty) of network state variables, such as pressures. Finally, it enables topology optimization to decouple network parts and prevent uncertainty propagation through the network. This is demonstrated using illustrative numerical experiments.
\end{abstract}


\thanksto{\noindent Vladimir Dvorkin and Dharik Mallapragada acknowledge funding from the MIT Energy Initiative Low-Carbon Energy Center on Electric Power Systems.}

\section{Introduction}

In many countries, balancing stochastic renewable generation in power systems relies on the flexibility of gas-fired power plants that couple power and natural gas systems. The disturbances from variable and uncertain renewables thus propagate to natural gas networks, causing  network  operators  to  adjust gas supply nominations to ensure   gas deliverability across the network. Such measures, however, are expensive and may not suffice to cope with weather-dependent gas extractions in the future. This motivates seeking alternative flexibility sources to control natural gas networks. 

Natural gas linepack -- temporally pressurized gas stored in pipelines -- has traditionally provided long-term flexibility for network operators (e.g., month-long storage), but more recently it has been modeled within shorter operational time frames. According to \cite{correa2014integrated} and \cite{schwele2019coordination}, linepack flexibility has a potential to save up to 1.4\% to 2.0\% of deterministic day-ahead scheduling cost. In the systems exposed to uncertainty, this cost saving potential increases up to 3.5\% on average \cite{ordoudis2019integrated}. The linepack flexibility also enables cost-security trade-offs in the stochastic dispatch of power and natural gas systems \cite{ratha2020affine}. 


The linepack models in \cite{ordoudis2019integrated} and \cite{ratha2020affine}, however, provide flexibility in the interest of stochastic day-ahead \textit{scheduling}. As two-stage stochastic programs, they co-optimize nominal set-points and reserve margins that immunize future operations against uncertainty. Unfortunately, those models do not answer how to \textit{control} the uncertain system state (e.g., nodal pressures and gas flows) as that uncertainty gradually realizes throughout the control (operating) horizon. Such a control task requires a \textit{multi-stage} stochastic modeling of network state, i.e., optimizing the network state as a function of all prior uncertainty realizations at any given stage, which is computationally challenging \cite{shapiro2005complexity}. This problem is thus typically limited to two subsequent stages \cite{dvorkin2020stochastic} or enjoys overly-conservative robust optimization techniques that disregard the temporal evolution of the stochastic renewable generation process \cite{roald2020uncertainty}.

In contrast to scheduling problems, this paper develops multi-stage stochastic control policies for natural gas networks with linepack and offers their computationally efficient optimization. Specifically, the following contributions are made:
\begin{enumerate}[labelwidth=!, labelindent=5pt]
    \item We develop multi-stage policies for controllable network assets that invoke linepack flexibility to guide operations within a finite discrete control horizon. These policies are based on multi-stage linear decision rules \cite{kuhn2011primal} that model uncertain network state through random trajectories and produce control inputs to maintain those trajectories within technical limits, hence expanding a two-stage stochastic control horizon modeled in \cite{dvorkin2020stochastic}. 
    \item We provide a distributionally robust chance-constrained optimization problem and its second-order cone program (SOCP) reformulation to optimize multi-stage control policies. We use the exact reformulation of double-sided distributionally robust constraints from \cite{xie2017distributionally} to reduce the conservatism of models in \cite{ratha2020affine,dvorkin2020stochastic}. 
    \item We provide several applications for the proposed multi-stage control policies in the sequential system coordination setting, where power system dispatch is followed by the gas network optimization problem. We find out that:
    \begin{itemize}
    \item[a)] The linepack flexibility suffices to balance stochastic extractions without substantial adjustments of the nominal gas injections. Invoking linepack flexibility in a 48-node network, we show that the standard deviation of stochastic gas extractions of 7.2\% of their nominal value can be accommodated while limiting the standard deviation of gas injections to 2.5\% of their nominal value. The cost-saving potential of linepack storage in our setting amounts to 10.3\% in expectation. 
    \item[b)] The proposed multi-stage control policies can be optimized to simultaneously minimize the variability (due to intermittency of renewable generation) and the variance (due to uncertainty) of the network state, thus improving on the variance-only control in \cite{dvorkin2020stochastic}. 
    \item[c)] Recognizing the role of network effects, the proposed control policies are additionally co-optimized with the binary valve deployment to identify the optimal network topology, which limits spatio-temporal disturbance propagations within a control horizon. 
    \end{itemize}
\end{enumerate}

\textit{Paper Organization:} Section \ref{sec:preliminaries} introduces the gas network equations and uncertainty model. Section  \ref{sec:Stoch_Opt} presents a chance-constrained multi-stage policy optimization problem, its tractable reformulation, and several applications. Section \ref{sec:experiments} provides numerical experiments and Section \ref{sec:conclusion} concludes. 

\textit{Notation:} lower- (upper-) case letters denote column vectors (matrices). $\mathbb{0}$ and $\mathbb{1}$ denote vectors (matrices) of zeros and ones; when ambiguous, their dimension is provided with a subscript. Operator $\text{dg}[x]$ returns a diagonal $n\times n$ matrix with diagonal entries of $n-$dimensional vector $x$, 
operator $[A]_{i}$ returns an $i$\textsuperscript{th} row $(1 \times n)$ of matrix $A$, operator $\text{Tr}[A]$ denotes matrix trace, and symbol $^\top$ stands for transposition. $\norm{\cdot}$ denotes $\ell_{2}-$norm. Boldface font denotes random variables and $\tilde{x}$ denotes the dependency of $x$ on random variables.   

\section{Preliminaries}\label{sec:preliminaries}

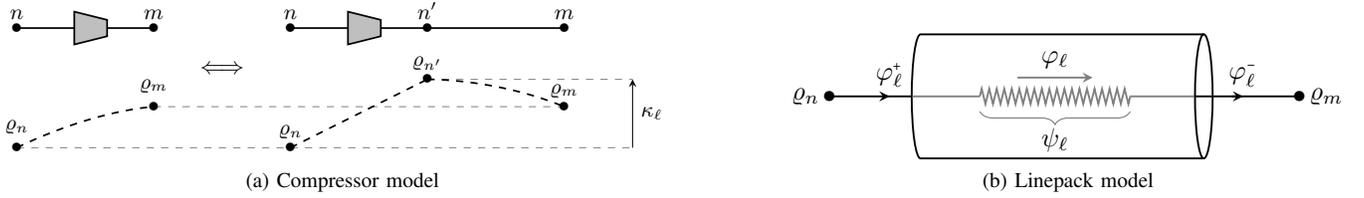
\begin{figure*}[t!] 
    \centering
    \subfloat[Compressor model\label{fig:compressor}]{%
    \resizebox{0.5\textwidth}{!}{%
    \begin{tikzpicture}
    \draw[gray,dashed] (0,-1.75) -- (4,-1.75);
    \draw[gray,dashed] (2,-1.15) -- (8,-1.15);
    \draw[gray,dashed] (4,-1.75) -- (9,-1.75);
    \draw[gray,dashed] (6,-0.75) -- (9,-0.75);
    \draw[line width=0.15mm, ->, >= stealth, black] (9,-1.75) -- (9,-0.75) node[midway,right]{\small $\kappa_{\ell}$};
    \draw[line width=0.25mm, black] (0,0) node[] {\small $\bullet$} node[above]{\small $n$} -- (2,0) node[] {\small $\bullet$} node[above]{\small $m$};
    \compressortwo{(0,0)}{0};
    \begin{scope}[xshift = 0cm, yshift = 0cm]
    \draw[dashed,line width=0.25mm] node[yshift=-1.75cm] {\small $\bullet$} node[yshift=-1.45cm] {\small $\varrho_{n}$} plot[domain=0:2,yshift=-0.75cm]   (\x,{0.5*\x-exp(0.1*\x^1.75)}) node[yshift=-0.75cm] {\small $\bullet$} node[yshift=-0.45cm] {\small $\varrho_{m}$};
    \end{scope}
    \draw[black] (3.0,-0.6) node {$\Longleftrightarrow$};
    \draw[line width=0.25mm, black] (4,0) node[] {\small $\bullet$} node[above]{\small $n$} -- (6,0) node[] {\small $\bullet$} node[above]{\small $n'$};
    \compressortwo{(4,0)}{0};
    \draw[line width=0.25mm, black] (6,0) -- (8,0) node[] {\small $\bullet$} node[above]{\small $m$};
    \begin{scope}[xshift = 4cm, yshift = -1.75cm]
    \draw[dashed,line width=0.25mm] node[above] {\small $\varrho_{n}$} node[] {\small $\bullet$} plot[domain=0:2]   (\x,{0.5*\x}) node[above] {\small $\varrho_{n'}$} node[] {\small $\bullet$};
    \end{scope}
    \begin{scope}[xshift = 6cm, yshift = 0.25cm]
    \draw[dashed,line width=0.25mm] plot[domain=0:2] (1*\x,{-exp(0.1*\x^1.75)}) node[above] {\small $\varrho_{m}$} node[] {\small $\bullet$};
    \end{scope}
    \end{tikzpicture}
    }
    }
    \hfill
    \subfloat[Linepack model \label{fig:linepack}]{%
    \begin{tikzpicture}
    \node[cylinder, line width=0.25mm, black, shape border rotate=0, draw,minimum height=4cm,minimum width=1.65cm]
    {};
    \draw[line width=0.25mm,gray] (-1.9,0)  -- (-1,0);
    \draw[line width=0.25mm,gray] (1,0)  -- (1.85,0);
    \draw[decorate,decoration={zigzag,segment length=1mm, amplitude=1mm},line width=0.25mm,gray] (-1,0) -- (1,0);
    \draw[line width=0.25mm,gray,->,>=stealth] (-0.5,0.25) -- node[above] {$\textcolor{black}{\varphi_{\ell}}$} (0.5,0.25);
    \draw [decorate,decoration={brace,amplitude=5pt},xshift=0pt,yshift=0pt,line width=0.15mm,gray]
    (1.0,-0.2) -- (-1.0,-0.2) node [midway,below,yshift = -0.1cm] {
    $\textcolor{black}{\psi_{\ell}}$};
    \draw[line width=0.25mm,black] (-2.25,0)  -- (-1.9,0);
    \draw[line width=0.25mm,black,->,>=stealth] (-3.0,0)  node[] {$\bullet$} node[left] {$\varrho_{n}$} -- (-2.20,0) node[above] {$\varphi_{\ell}^{\shortplus}$};
    \draw[line width=0.25mm,black] (2.35,0) -- (3.25,0) node[] {$\bullet$} node[right] {$\varrho_{m}$};
    \draw[line width=0.25mm,black,->,>=stealth] (1.85,0)  -- (2.5,0) node[above] {$\varphi_{\ell}^{\shortminus}$};
    \end{tikzpicture}
    }
    \caption{Modeling of pressure regulation (left) and linepack natural gas storage (right) and the corresponding notation.}
    \vspace{-0.6cm}
\end{figure*}    

\subsection{Natural Gas Network Equations}\label{subsec:gas_equations}
Consider a gas network with the sets $\mathcal{N}=\{1,\dots,N\}$ of nodes and $\mathcal{E}=\{1,\dots,E\}$ of edges (pipelines) connecting those nodes. Each edge is assigned a direction from sending end $n$ to receiving end $m$, i.e., if $(n,m)\!\in\!\mathcal{E}$, then $(m,n)\!\notin\!\mathcal{E}$, but the flow is allowed to be negative. This topology is encoded in the node-edge incidence matrix $A\!\in\!\mathbb{R}^{N\times E}$ as
\begin{align*}
    A_{k\ell} = 
    \left\{
    \begin{array}{rl}
        +1, & \text{if}\;k=n  \\
        -1, & \text{if}\;k=m \\
        0, & \text{otherwise}
    \end{array}
    \right.\quad\forall\ell=(n,m)\in\mathcal{E}.
\end{align*}
The gas network optimization problem identifies optimal gas injections $\vartheta\in\mathbb{R}^{N}$ that satisfy gas extractions $\delta\in\mathbb{R}_{+}^{N}$ across the network. Towards this goal, the network operator must maintain such nodal pressures $\varrho\in\mathbb{R}^{N}$ that do not prevent the gas flow $\varphi\in\mathbb{R}^{E}$ in the network. To support pressures, the network operator deploys active pipelines $\mathcal{E}_{a}\!\subseteq\!\mathcal{E}$ that host compressors $\mathcal{E}_{c}\!\subseteq\!\mathcal{E}_{a}$ or control valves $\mathcal{E}_{v}\!\subseteq\!\mathcal{E}_{a}$, $\mathcal{E}_{c}\cap\mathcal{E}_{v}=\emptyset$, that support pressure in a continues manner \cite{fugenschuh2015chapter}. A compressor edge $\ell=(n,m)\!\in\!\mathcal{E}_{c}$ is modeled as in Fig. \ref{fig:compressor}, where the controllable compression rate $\kappa_{\ell}\geqslant0$ is followed by the non-linear pressure drop according to the Weymouth equation \cite{fugenschuh2015chapter}. A valve edge $\ell\in\mathcal{E}_{v}$ is modeled similarly, but the linear decompression rate $\kappa_{\ell}$ always remains non-positive. 
To regulate pressure, active pipelines consume gas according to conversion factors stored in matrix $B\in\mathbb{R}^{N\times E}$:
\begin{align*}
    B_{k\ell} = 
    \left\{
    \begin{array}{rl}
        b_{\ell}, & \text{if}\;k=n,\;k\in\mathcal{E}_{c}  \\
        -b_{\ell}, & \text{if}\;k=m,\;k\in\mathcal{E}_{v}  \\
        0, & \text{otherwise}
    \end{array}
    \right.\quad\forall\ell=(n,m)\in\mathcal{E},
\end{align*}
where $b_{\ell}$ is a conversion factor from the regulation rate to gas mass. This gives rise to the following gas flow equations
\begin{subequations}\label{eq:network_equations}
\begin{align}
    &A\varphi = \vartheta - B\kappa - \delta,\label{eq:conserv_law}\\
    &\varphi_{\ell}|\varphi_{\ell}|=w_{\ell}\braceit{(\varrho_{n} + \kappa_{\ell})^2 - \varrho_{m}^2},&&\forall\ell=(n,m)\in\mathcal{E}, \label{eq:weymouth_eq}\\
    &\varphi_{\ell}\geqslant 0,\;\forall\ell\in\mathcal{E}_{a},\label{eq:uni_direction}
\end{align}
where equation \eqref{eq:conserv_law} ensures gas mass conservation by balancing gas flows, injections and extractions, equation \eqref{eq:weymouth_eq} models the gas flow in pipelines using the Weymouth equation, where $w\in\mathbb{R}_{+}^{E}$ is the friction coefficient, and \eqref{eq:uni_direction} models unidirectional gas flow required for active pipelines only. 

To model linepack, we adopt the notation in Fig. \ref{fig:linepack}, where $\varphi_{\ell}^{\shortplus},\varphi_{\ell}$ and $\varphi_{\ell}^{\shortminus}$ respectively denote the inlet, midway and outlet gas flows, and $\psi_{\ell}$ denotes the amount of linepack. The linepack model in \cite{schwele2019coordination} couples these variables as
\begin{align}
    &\varphi_{\ell} = \tfrac{1}{2}\braceit{\varphi_{\ell}^{\shortplus} + \varphi_{\ell}^{\shortminus}},\label{eq:mean_flow}\\
    &\psi_{\ell} = \tfrac{1}{2}s_{\ell}\braceit{\varrho_{n}+\kappa_{\ell}+\varrho_{m}},\label{eq:line_pack_def}
\end{align}
where the first entry sets the midway flow as an average flow through the pipeline, and the second entry defines the linepack proportionally to the pressures at the sending and receiving ends, with $s_{\ell}\geqslant0$ being a pipe-specific conversion factor. Observe, the linepack equations \eqref{eq:mean_flow} and \eqref{eq:line_pack_def} are related to the Weymouth equation \eqref{eq:weymouth_eq}, which requires the midway flow to be proportional to the loss of pressure. The linepack flexibility calls for the following relaxation of conservation law \eqref{eq:conserv_law}:
\begin{align}
    &A^{\shortplus}\varphi^{\shortplus} + A^{\shortminus}\varphi^{\shortminus}
    =\vartheta-B\kappa-\delta,\label{det:conservation_law_lp}
\end{align}
which requires flow conservation only in terms of inlet and outlet gas flows, where matrices $A^{\shortplus}$ and $A^{\shortminus}$ are obtained from the full incidence matrix $A$ by striking out only negative and positive entries, respectively. To model linepack dynamics across a finite set $\mathcal{T}=\{1,\dots,T\}$ of time stages, we consider
\begin{align}
    &\psi_{t} = \psi_{t-1} + \varphi_{t}^{\shortplus} - \varphi_{t}^{\shortminus}, \quad\forall t\in\mathcal{T},\label{eq:lp_in_time}\\
    &\psi_{T}\geqslant\psi_{\ell0}, \label{eq:non_depletion}
\end{align}
where equation \eqref{eq:lp_in_time} models the linepack state in time and equation \eqref{eq:non_depletion} does not permit the linepack depletion at the last stage $(t=T)$ below the initial amount $(t=0)$. 
\end{subequations}

\subsection{Linearization of the Weymouth Equation \eqref{eq:weymouth_eq}}
When gas extractions are uncertain, solving the system of equations \eqref{eq:network_equations} is difficult due to the non-convex Weymouth equation. We thus provide the following linearization technique to reduce problem complexity.
Let $\mathcal{J}(x)\in\mathbb{R}^{E\times m}$ denote the Jacobian of \eqref{eq:weymouth_eq} in some arbitrary point $x\in\mathbb{R}^{m}$, and let $(\mathring{\varphi},\mathring{\varrho},\mathring{\kappa})$ be a stationary point, i.e., obtained by solving equations \eqref{eq:network_equations} for the deterministic (mean) gas extraction values. Then, according to \cite{dvorkin2020stochastic}, the following linear equation
\begin{align}
    \varphi &= 
    \underbrace{\mathcal{J}(\mathring{\varphi})^{-1}(\mathcal{J}(\mathring{\varrho})\mathring{\varrho} + \mathcal{J}(\mathring{\kappa})\mathring{\kappa}) + \mathring{\varphi}}_{w_{0}(\mathring{\varphi},\mathring{\varrho},\mathring{\kappa})} \nonumber \\
    &\quad\quad\quad\underbrace{-\mathcal{J}(\mathring{\varphi})^{-1}\mathcal{J}(\mathring{\varrho})}_{W_{1}(\mathring{\varphi},\mathring{\varrho})}\varrho
    \underbrace{-\mathcal{J}(\mathring{\varphi})^{-1}\mathcal{J}(\mathring{\kappa})}_{W_{2}(\mathring{\varphi},\mathring{\kappa})}\kappa\nonumber\\
    &=w_{0} + W_{1}\varrho + W_{2}\kappa\label{eq:wey_lin}
\end{align}
is equivalent to equation \eqref{eq:weymouth_eq} at the stationary point, where $w_{0}\in\mathbb{R}^{E}$, $W_{1}\in\mathbb{R}^{E\times N}$ and $W_{2}\in\mathbb{R}^{E\times E}$ are sensitivities depending on the stationarity point. To ensure the flow solution uniqueness under linearization \eqref{eq:wey_lin}, we require the pressure at the reference node $\text{r}$ (e.g., with a large constant gas injection) to be fixed as $\varrho_{\text{r}t}=\mathring{\varrho}_{\text{r}t}$. When deviations from the stationary point are large, however, linearization \eqref{eq:wey_lin} may produce approximation errors. The work in \cite{dvorkin2020stochastic}, however, bounds these errors to modest values even under large deviations. 

\subsection{From Renewable Power to Gas Extraction Uncertainty}\label{subsec:from_wind_to_gas}

The gas extraction uncertainty is due to combined cycle gas turbines (CCGT) that balance stochastic renewable generation in power systems. To model renewable uncertainty, consider a dynamic stochastic process $\fat{\zeta}\!=\!(\fat{\zeta}_{1},\dots,\fat{\zeta}_{T})\in\mathbb{R}^{k}$  which develops throughout a $T-$stage horizon. Here, vector $\fat{\zeta}_{t}\in\mathbb{R}^{k_{t}}$ collects random variables revealed at stage $t$. Let vector $\boldsymbol{\zeta}^{t}=(\boldsymbol{\zeta}_{1},\dots,\boldsymbol{\zeta}_{t})\in\mathbb{R}^{k^{t}}$ collect random variables realized upon stage $t$, i.e., current and all previous realizations, meaning $k^{t}=\sum_{\tau=1}^{t}k_{\tau}$ and $k^{T}=k$. The random generation from a set $\mathcal{R}=\{1,\dots,R\}$ of renewable units is modeled as
\begin{align}\label{eq:renewable_process}
    \tilde{r}_{t}(\fat{\zeta}^{t}) = \Omega_{t}\fat{\zeta}^{t} \in\mathbb{R}^{R},
\end{align}
i.e., as an affine combination of random variables revealed up to stage $t$, where $\Omega_{t}\in\mathbb{R}^{R\times k^{t}}$ is a matrix of coefficients. We consider that  $\fat{\zeta}$ belongs to the unknown distribution $\mathbb{P}_{\fat{\zeta}}(\widehat{\mu},\widehat{\Sigma})$, but its mean $\widehat{\mu}\in\mathbb{R}^{k}$ and covariance $\widehat{\Sigma}\in\mathbb{R}^{k\times k}$ are known from history. Finally, we set $\fat{\zeta}_{1}\!=\!1$ $(k_{1}\!=\!1)$ as the renewable generation at time stage 1 (stage ``here-and-now'') is certain. 

To map renewable uncertainty $\fat{\zeta}$ into uncertain gas extractions, suppose that the power system operator solves a multi-stage optimal power flow problem, where the affine generator response $\tilde{g}_{t}(\fat{\zeta}^{t}) = G_{t}\fat{\zeta}^{t}\in\mathbb{R}^{M}$ is computed for all $M$ generators by optimizing variable $G_{t}\in\mathbb{R}^{M\times k^{t}}$, which is consistent with \cite{bienstock2014chance}. Then, for a power-to-gas conversion matrix $\Lambda\in\mathbb{R}^{M\times N}$, the uncertain extraction is modeled as $\tilde{\delta}_{t}(\fat{\zeta}^{t})=\Lambda G_{t}^{\star} \fat{\zeta}^{t}\in\mathbb{R}^{N}$, i.e., as an affine combination of random renewable deviations, where $G_{t}^{\star}$ is an optimized generator response. To ease the narrative, we make a substitution $\Delta_{t} = \Lambda G_{t}^{\star}\in\mathbb{R}^{N\times k^{t}}$. Once the generator response is optimized, we assume that the power system operator submits uncertainty information -- matrices $\{\Delta_{t}\}_{t\in\mathcal{T}}$ and two moments $\widehat{\mu}$ and $\widehat{\Sigma}$ -- to the gas network operator.  

\begin{remark*}\normalfont 
In the interest of space, we do not discuss the optimization of generator policy $\tilde{g}_{t}(\fat{\zeta}^{t})$ explicitly, but provide it in the online repository \cite{dvorkin2021}.  
\end{remark*}

\section{Invoking Linepack Flexibility through Multi-Stage Stochastic Optimization}\label{sec:Stoch_Opt}

The uncertainty of gas extractions motivates the following stochastic formulation of the gas network optimization: 
\begin{subequations}\label{prob:stoch_intr}
\begin{align}
&\text{minimize}_{\tilde{\vartheta}_{t},\tilde{\kappa}_{t},\tilde{\varphi}_{t},\tilde{\varrho}_{t},\tilde{\varphi}_{t}^{\shortplus},\tilde{\varphi}_{t}^{\shortminus},\tilde{\psi}_{t}}\nonumber\\
&\mathbb{E}_{\mathbb{P}_{\fat{\zeta}}}\!\!\left[\sum_{t=1}^{T}\braceit{c_{1}^{\top}\tilde{\vartheta}_{t}(\fat{\zeta}^{t}) + \tilde{\vartheta}_{t}(\fat{\zeta}^{t})^{\top}\text{dg}[c_{2}]\tilde{\vartheta}_{t}(\fat{\zeta}^{t})}\right]\label{stoch_intr_obj}\\
&\text{subject to}\colon\nonumber\\
&\mathbb{P}_{\fat{\zeta}^{t}}\!\!\left[
\begin{aligned}
&A^{\shortplus}\tilde{\varphi}_{t}^{\shortplus}(\fat{\zeta}^{t}) + A^{\shortminus}\tilde{\varphi}_{t}^{\shortminus}(\fat{\zeta}^{t})
=\tilde{\vartheta}_{t}(\fat{\zeta}^{t})\\
&\qquad\qquad\qquad\qquad\;\!-B\tilde{\kappa}_{t}(\fat{\zeta}^{t})-\tilde{\delta}_{t}(\fat{\zeta}^{t}) \\
&\tilde{\varphi}_{t}(\fat{\zeta}^{t})=w_{0t} + W_{1t}\tilde{\varrho}_{t}(\fat{\zeta}^{t}) + W_{2t}\tilde{\kappa}_{t}(\fat{\zeta}^{t})\\
&\tilde{\varrho}_{\text{r}t}(\fat{\zeta}^{t})=\mathring{\varrho}_{\text{r}t}\\
&\tilde{\varphi}_{t}(\fat{\zeta}^{t}) = \tfrac{1}{2}\braceit{\tilde{\varphi}_{t}^{\shortplus}(\fat{\zeta}^{t}) + \tilde{\varphi}_{t}^{\shortminus}(\fat{\zeta}^{t})}\\
&\tilde{\psi}_{t}(\fat{\zeta}^{t}) = \tfrac{1}{2}\text{dg}[s]\braceit{\tilde{\kappa}_{t}(\fat{\zeta}^{t}) + |A|^{\top}\tilde{\varrho}_{t}(\fat{\zeta}^{t})}\\
&\tilde{\psi}_{t}(\fat{\zeta}^{t}) = \tilde{\psi}_{t\shortminus1}(\fat{\zeta}^{t\shortminus1}) + \tilde{\varphi}_{t}^{\shortplus}(\fat{\zeta}^{t}) - \tilde{\varphi}_{t}^{\shortminus}(\fat{\zeta}^{t})
\end{aligned}
\right]=1,\label{stoch_intr_as}
\\
&\mathbb{P}_{\fat{\zeta}^{t}}\!\!\left[
\begin{aligned}
\begin{aligned}
&\underline{\vartheta}\leqslant\tilde{\vartheta}_{t}(\fat{\zeta}^{t})\leqslant\overline{\vartheta}\\
&\underline{\varrho}\leqslant \tilde{\varrho}_{t}(\fat{\zeta}^{t})\leqslant\overline{\varrho}\\
&\underline{\kappa}\leqslant \tilde{\kappa}_{t}(\fat{\zeta}^{t})\leqslant\overline{\kappa}\\
\end{aligned} 
& &
\begin{aligned}
&\tilde{\psi}_{T}(\fat{\zeta}^{T})\geqslant\psi_{0}\\
&\tilde{\varphi}_{t\ell}(\fat{\zeta}^{t})\geqslant0,\;\forall\ell\in\mathcal{E}_{a}
\end{aligned} 
\end{aligned}
\right]\!\geqslant\!1\!-\!\varepsilon_{t},\label{stoch_intr_cc}\\
&\text{for all time stages}\;t\in\mathcal{T}.\nonumber
\end{align}
\end{subequations}
This problem minimizes the expected quadratic cost of gas injection, with cost vectors $c_{1}$ and $c_{2}$, subject to probabilistic constraints enforced across the $T-$stage control horizon. Here, the stochastic gas injection $\tilde{\vartheta}_{t}$ and pressure regulation rate $\tilde{\kappa}_{t}$ are variable vectors controlled by the network operator, and stochastic pressure $\tilde{\varrho}_{t}$, flow $\tilde{\varphi}_{t}^{(\cdot)}$ and linepack $\tilde{\psi}_{t}$ vectors are state variables that adjust in response to the control variables. Since the pressure at the reference node is fixed to the reference value $\mathring{\varrho}_{\text{r}t}$, the response of state variables is unique \cite{dvorkin2020stochastic}.  Constraint \eqref{stoch_intr_as} requires the stochastic equalities to hold with probability 1, i.e., for any realizations of uncertainty. The entries of this constraint are stochastic counterparts of the gas network equations introduced in Section \ref{sec:preliminaries}. $T$ joint chance constraints in \eqref{stoch_intr_cc} are introduced to ensure the satisfaction of minimal and maximum physical limits with probability at least $1-\varepsilon_{t}$, where $\varepsilon_{t}$ is a small prescribed parameter.

Unlike two-stage scheduling problems under uncertainty, the variables and constraints in formulation \eqref{prob:stoch_intr} are required to depend on $\fat{\zeta}^{t}=(\boldsymbol{\zeta}_{1},\dots,\boldsymbol{\zeta}_{t})$ -- all prior realization of renewable forecast errors up to stage $t$. This enables expressing the network state through stochastic \textit{trajectories} that must be controlled within network limits up to parameter $\varepsilon_{t}$ at every time stage $t$. This formulation, however, increases computational complexity (exponentially rather than linearly) as $T$ increases \cite{shapiro2005complexity}. Fortunately,  this multi-stage stochastic program can be approximated efficiently under linear decision rules \cite{kuhn2011primal}.  

\subsection{Reformulations via Linear Decision Rules}
To reformulate the multi-stage stochastic problem \eqref{prob:stoch_intr}, we introduce truncation operators $S_{t},\forall t\in\mathcal{T}$, defined in \cite{kuhn2011primal} as
\begin{align*}
    S_{t}:\mathbb{R}^{k}\mapsto\mathbb{R}^{k^{t}},\quad \fat{\zeta}\mapsto\fat{\zeta}^{t},
\end{align*}
to extract only those random entries of vector $\fat{\zeta}$ that realize by a particular time stage $t$.
Using truncation operators, the variables of stochastic problem \eqref{prob:stoch_intr} are defined through the following multi-stage linear decision rules:
\begin{align}\label{eq:LDRs}
    \begin{aligned}
    &\begin{aligned}
    &\tilde{\vartheta}_{t}(\fat{\zeta}^{t})=\Theta_{t}S_{t}\fat{\zeta}\\
    &\tilde{\kappa}_{t}(\fat{\zeta}^{t})=K_{t}S_{t}\fat{\zeta}\\
    \end{aligned} 
    &\begin{aligned}
    &\tilde{\varrho}_{t}(\fat{\zeta}^{t})=P_{t}S_{t}\fat{\zeta}\\
    &\tilde{\psi}_{t}(\fat{\zeta}^{t})= \Psi_{t}S_{t}\fat{\zeta}\\
    \end{aligned}
    &
    &\begin{aligned}
    &\tilde{\varphi}_{t}(\fat{\zeta}^{t})=\Phi_{t}S_{t}\fat{\zeta}\\
    &\tilde{\varphi}_{t}^{\shortplus}(\fat{\zeta}^{t})=\Phi_{t}^{\shortplus}S_{t}\fat{\zeta}\\
    &\tilde{\varphi}_{t}^{\shortminus}(\fat{\zeta}^{t})=\Phi_{t}^{\shortminus}S_{t}\fat{\zeta}\\
    \end{aligned}
    \end{aligned},
\end{align}
where $\Theta_{t},P_{t}\!\in\!\mathbb{R}^{N\times k^{t}}$, $K_{t},\Psi_{t},\Phi_{t}^{(\cdot)}\!\in\!\mathbb{R}^{E\times k^{t}}$ are coefficient matrices subject to optimization. These matrices share an important property: their first column corresponds to the nominal value, e.g., gas supply nomination at stage $t$, and the remaining columns define the variable recourse with respect to all prior realizations of uncertainty up to stage $t$. This variable representation makes possible the following reformulations. 

\subsubsection{Objective Function Reformulation}\label{subsubsec:ref_cost} The expected cost in \eqref{stoch_intr_obj} is reformulated under linear decision rules \eqref{eq:LDRs} using truncation and the first two moments of distribution $\mathbb{P}_{\fat{\zeta}}$ as 
\begin{align}
    &\mathbb{E}_{\mathbb{P}_{\fat{\zeta}}}\left[\sum_{t=1}^{T}\braceit{c_{1}^{\top}\Theta_{t}S_{t}\fat{\zeta} + \braceit{\Theta_{t}S_{t}\fat{\zeta}}^{\top}\text{dg}[c_{2}]\Theta_{t}S_{t}\fat{\zeta}}\right]\Longrightarrow\nonumber\\
    &\sum_{t=1}^{T}\!\braceit{c_{1}^{\top}\Theta_{t}S_{t}\widehat{\mu} + \text{Tr}\left[
    \Theta_{t}^{\top}\text{dg}[c_{2}]\Theta_{t}S_{t}\braceit{\widehat{\Sigma} + \widehat{\mu}\widehat{\mu}^{\top}}S_{t}^{\top}
    \right]},
\end{align}
which is a convex expression easy to optimize. Here, the second term is obtained knowing that the expected value of the outer product of $\fat{\zeta}$ amounts to $\mathbb{E}_{\mathbb{P}_{\fat{\zeta}}}[\fat{\zeta}\fat{\zeta}^{\top}]=\widehat{\Sigma}+\widehat{\mu}\widehat{\mu}^{\top}$. 

\subsubsection{Reformulation of Stochastic Equalities}\label{subsubsec:ref_as} The entries of constraint \eqref{stoch_intr_as} admit an analytic reformulation by constraining the matrices in \eqref{eq:LDRs}: each entry can be represented as a linear stochastic equality $X_{t}S_{t}\fat{\zeta}=\mathbb{0}_{m}$ for some coefficient matrix $X_{t}\in\mathbb{R}^{m\times k^{t}}$. One then verifies that this equality holds for any $\fat{\zeta}$ if $X_{t}$ is constrained by $X_{t}S_{t}=\mathbb{0}_{m\times k}$. For instance, the pressure at the reference node reformulates accordingly as
\begin{align*}
    &\tilde{\varrho}_{\text{r}t}(\fat{\zeta}^{t})=\mathring{\varrho}_{\text{r}t} 
    \overset{\fat{\zeta}_{1}=1}{\Longleftrightarrow}
    [P_{t}S_{t}\fat{\zeta}]_{\text{r}}=\begin{bmatrix}\mathring{\varrho}_{\text{r}t}&\mathbb{0}_{1\times (k^{t}\shortminus1)}\end{bmatrix} S_{t}\fat{\zeta}
    \\
    &\Longleftrightarrow\braceit{[P_{t}]_{\text{r}} - \begin{bmatrix}\mathring{\varrho}_{\text{r}t}&\mathbb{0}_{1\times (k^{t}\shortminus1)}\end{bmatrix}}S_{t}=\mathbb{0}.
\end{align*}
The remaining stochastic equalities are reformulated by analogy, whose final expressions are given in \eqref{prob_cc_eq_1}-\eqref{prob_cc_eq_6}. 

\subsubsection{Reformulation of Joint Chance Constraints} \label{subsubsec:ref_cc} We follow the analytic reformulation based on Bonferroni and Chebyshev approximations. That is, we split $T$ joint chance constraints \eqref{stoch_intr_cc}  into $(E + |\mathcal{E}_{a}|) + (N+2E)$ individual chance constraints each, where the first cardinality term relates to single-sided constraints and the second one to double-sided constraints. Single-sided constraints are treated using a distributionally robust Chebyshev reformulation \cite{xie2017distributionally}. For instance, the constraints on the minimum stochastic gas flow in active pipelines
\begin{align*}
    &\mathbb{P}_{\fat{\zeta}^{t}}\left[
    \tilde{\varphi}_{t\ell}(\fat{\zeta}^{t})\geqslant0
    \right]\geqslant1-\overline{\varepsilon}_{t}, \;\forall \ell\in\mathcal{E}_{a}, \forall t\in\mathcal{T},
\end{align*}
are reformulated into second-order cone constraint \eqref{prob_cc_soc_min_flow}, where $\overline{\varepsilon}_{t} = (1-\varepsilon_{t})/(E + |\mathcal{E}_{a}| + N+2E)$ and $\widehat{F}$ is covariance matrix factorization, i.e., $\widehat{\Sigma}=\widehat{F}\widehat{F}^{\top}$. Notice, a similar reformulation of double-sided chance-constraints, e.g., on stochastic gas injections bounded from both ends, is economically inefficient because the lower and upper bounds cannot be violated simultaneously. We thus invoke the exact reformulation of such constraints from \cite[Th. 2]{xie2017distributionally}. For example, the chance constraint on stochastic gas injection
\begin{align*}
    &\mathbb{P}_{\fat{\zeta}^{t}}\left[
    \underline{\vartheta}_{n}\leqslant[\Theta_{t}]_{n}S_{t}\fat{\zeta}\leqslant\overline{\vartheta}_{n}
    \right]\geqslant1-\overline{\varepsilon}_{t}, \;\forall n\in\mathcal{N},\;\forall t\in\mathcal{T},
\end{align*}
re-writes into a second-order cone constraint \eqref{prob_cc_gas_inj_soc} and a set of linear constraints in \eqref{prob_cc_gas_inj_lin_1} and \eqref{prob_cc_gas_inj_lin_2}, where $y_{tn}^{\vartheta}$ and $x_{tn}^{\vartheta}$ are auxiliary variables subject to optimization. The remaining double-sided chance-constraints are reformulated by analogy.

The tractable reformulation of the multi-stage stochastic problem \eqref{prob:stoch_intr} in linear decision rules \eqref{eq:LDRs} takes the form:
\begin{subequations}\label{prob:stoch_tr}
\begin{align}
&\text{minimize}_{\Theta_{t},P_{t},K_{t},\Psi_{t},\Phi_{t}^{(\cdot)},x_{t}^{(\cdot)},y_{t}^{^{(\cdot)}}}\nonumber\\
    &\sum_{t=1}^{T}\!\braceit{c_{1}^{\top}\Theta_{t}S_{t}\widehat{\mu}\!+\!\text{Tr}\left[
    \Theta_{t}^{\top}\text{dg}[c_{2}]\Theta_{t}S_{t}\braceit{\widehat{\Sigma}\!+\!\widehat{\mu}\widehat{\mu}^{\top}}S_{t}^{\top}
    \right]}\label{prob_cc_obj}\\
    &\text{subject to}\colon\nonumber\\
    &\braceit{A^{\shortplus}\Phi_{t}^{\shortplus} + A^{\shortminus}\Phi_{t}^{\shortminus} - \Theta_{t}S_{t}+BK_{t}+\Delta_{t}}S_{t} = \mathbb{0}\label{prob_cc_eq_1}\\
    &\braceit{\Phi_{t} - \big[w_{0t}\;\;\mathbb{0}_{E\times (k^{t}\shortminus1)}\big] -  W_{1t}P_{t} - W_{2t}K_{t}}S_{t} = \mathbb{0} \\
    &\braceit{[P_{t}]_{\text{r}} - \big[\mathring{\varrho}_{\text{r}t}\;\;\mathbb{0}_{1\times (k^{t}\shortminus1)}\big]}S_{t}=\mathbb{0}\\
    &\braceit{\Phi_{t} - \tfrac{1}{2}\Phi_{t}^{\shortplus} - \tfrac{1}{2}\Phi_{t}^{\shortminus}}S_{t}=\mathbb{0}\\
    &\braceit{\Psi_{t}-\tfrac{1}{2}\text{dg}[s]\braceit{K_{t} + |A|^{\top}P_{t}}}S_{t}=\mathbb{0}\\
    &\braceit{\Psi_{t}S_{t} - \Phi_{t}^{\shortplus} + \Phi_{t}^{\shortminus}}S_{t} - \Psi_{(t\shortminus1)}S_{(t\shortminus1)}=\mathbb{0}\label{prob_cc_eq_6}\\
    &\textstyle\sqrt{\frac{1-\overline{\varepsilon}_{t}}{\overline{\varepsilon}_{t}}}\norm{
    \widehat{F}
    [\Phi_{t}S_{t}]_{\ell}^{\top}
    }\leqslant
    [\Phi_{t}S_{t}\widehat{\mu}]_{\ell}\;^* \label{prob_cc_soc_min_flow}\\
    &\textstyle\sqrt{\frac{1-\overline{\varepsilon}_{t}}{\overline{\varepsilon}_{t}}}\norm{
    \widehat{F}
    [\Psi_{t}S_{t}]_{\ell}^{\top}
    }\leqslant
    [\Psi_{t}S_{t}\widehat{\mu} - \psi_{0}]_{\ell} \\
    &\sqrt[-2]{\overline{\varepsilon}}\norm{\begin{matrix}\widehat{F}[\Theta_{t}S_{t}]_{n}^{\top}\\y_{tn}^{\vartheta}\end{matrix}}\leqslant\textstyle\frac{1}{2}\big(\overline{\vartheta}_{n} - \underline{\vartheta}_{n}\big) - x_{tn}^{\vartheta}\label{prob_cc_gas_inj_soc}\\
    &\sqrt[-2]{\overline{\varepsilon}}\norm{\begin{matrix}\widehat{F}[K_{t}S_{t}]_{\ell}^{\top}\\y_{t\ell}^{\kappa}\end{matrix}}\leqslant\textstyle\frac{1}{2}\big(\overline{\kappa}_{\ell} - \underline{\kappa}_{\ell}\big) - x_{t\ell}^{\kappa}\\
    &\sqrt[-2]{\overline{\varepsilon}}\norm{\begin{matrix}\widehat{F}[P_{t}S_{t}]_{n}^{\top}\\y_{tn}^{\varrho}\end{matrix}}\leqslant\textstyle\frac{1}{2}\big(\overline{\varrho}_{n} - \underline{\varrho}_{n}\big) - x_{tn}^{\varrho}\\
    &\left|[\Theta_{t}S_{t}]_{n}\widehat{\mu} - \textstyle\frac{1}{2}\big(\overline{\vartheta}_{n} - \underline{\vartheta}_{n}\big) \right| \leqslant y_{tn}^{\vartheta} + x_{tn}^{\vartheta}\label{prob_cc_gas_inj_lin_1}\\
    &\left|[K_{t}S_{t}]_{\ell}\widehat{\mu} - \textstyle\frac{1}{2}\big(\overline{\kappa}_{\ell} - \underline{\kappa}_{\ell}\big) \right| \leqslant y_{t\ell}^{\kappa} + x_{t\ell}^{\kappa}\\
    &\left|[P_{t}S_{t}]_{n}\widehat{\mu} - \textstyle\frac{1}{2}\big(\overline{\varrho}_{n} - \underline{\varrho}_{n}\big)\right| \leqslant y_{tn}^{\varrho} + x_{tn}^{\varrho}\\
    &\textstyle\frac{1}{2}\big(\overline{\vartheta}_{n} - \underline{\vartheta}_{n}\big)\geqslant x_{tn}^{\vartheta} \geqslant 0,\; y_{tn}^{\vartheta} \geqslant 0\label{prob_cc_gas_inj_lin_2}\\
    &\textstyle\frac{1}{2}\big(\overline{\kappa}_{\ell} - \underline{\kappa}_{\ell}\big)\geqslant x_{t\ell}^{\kappa} \geqslant 0,\; y_{t\ell}^{\kappa} \geqslant 0\\
    &\textstyle\frac{1}{2}\big(\overline{\varrho}_{n} - \underline{\varrho}_{n}\big)\geqslant x_{tn}^{\varrho} \geqslant 0,\; y_{tn}^{\varrho} \geqslant 0 \label{prob_cc_soc_last}\\
    &\forall t\in\mathcal{T},\;\forall n\in\mathcal{N},\;\forall\ell\in\mathcal{E}, \;^*\forall\ell\in\mathcal{E}_{a}.  
\end{align}
\end{subequations}

\subsection{Stochastic Control Applications}\label{subsec:control_applications}
The gas network operator controls the network state by optimizing stochastic gas injections $\tilde{\vartheta}_{t}(\fat{\zeta}^{t})$ and pressure regulation rates $\tilde{\kappa}_{t}(\fat{\zeta}^{t})$ across the $T-$stage control horizon. This subsection discusses how these control variables can be optimized by chance-constrained program \eqref{prob:stoch_tr} in several applications. 

\subsubsection{Linepack to Balance Stochastic Renewables}\label{subsubsec:linepack} Consider a network operator who accommodates uncertain gas extractions in real-time and is reluctant to adjust nominal gas injection contracts, because substantially altering gas nominations leads to expensive regulation penalties. Instead of substantially altering nominal injection, the network operator can leverage the linepack as an alternative source of flexibility to balance gas extractions by optimizing pressure regulation rates of compressors and valves. Doing so means optimizing program \eqref{prob:stoch_tr} while additionally ensuring the satisfaction of the following constraint on stochastic gas injection variables: 
\begin{subequations}
\begin{align}
    &\norm{\widehat{F}[\Theta_{t}]_{n}S_{t}}
    \leqslant \alpha^{\vartheta}[\Theta_{t}]_{n}S_{t}\widehat{\mu},\; \forall n\in\mathcal{N},
\end{align}
\end{subequations}
where the left-hand side computes the standard deviation $\text{Std}[\tilde{\vartheta}_{t}(\fat{\zeta}^{t})]$ of stochastic injections, and the right-hand side is the fraction of the nominal gas injection up to factor $\alpha^{\vartheta}\geqslant0$. Setting $\alpha^{\vartheta}=0.025$, for example, constrains the standard deviation of gas injections under the optimized control policies to 2.5\% of the nominal value, which is relatively small.



\subsubsection{Minimal Variability of Network State}\label{subsubsec:var}

Linepack flexibility can be invoked to minimize the variability of state variables, such as nodal gas pressures. Towards the goal, objective function \eqref{prob_cc_obj} is augmented with
\begin{subequations}
\begin{align}
  \mathbb{E}_{\mathbb{P}_{\fat{\zeta}}}\left[\alpha^{\varrho}\sum_{t=2}^{T}
  \big\lVert
  P_{t}S_{t}\fat{\zeta} - P_{t-1}S_{t-1}\fat{\zeta}
  \big\rVert
  \right],\label{eq:rho_var}
\end{align}
which minimizes the expected distance between stochastic nodal pressures at adjacent stages of control horizon up to prescribed non-negative penalty parameter $\alpha^{\varrho}$. Hence, by setting $\alpha^{\varrho}$, the control variables of \eqref{prob:stoch_tr} are optimized to meet a trade-off between the expected costs and variability of system state. To reformulate \eqref{eq:rho_var}, we make a substitution $\tilde{P}_{t}=P_{t}S_{t}$ for compactness. Then, this term becomes
\begin{align}
  &\mathbb{E}_{\mathbb{P}_{\fat{\zeta}}}\left[
  \alpha^{\varrho}\sum_{t=2}^{T}
  \big\lVert
  \tilde{P}_{t}\fat{\zeta} - \tilde{P}_{t-1}\fat{\zeta}
  \big\rVert
  \right] \nonumber\\
  &= \mathbb{E}_{\mathbb{P}_{\fat{\zeta}}}\left[
  \alpha^{\varrho}\sum_{t=2}^{T}
  \braceit{\braceit{\tilde{P}_{t} - \tilde{P}_{t-1}}\fat{\zeta}}^{\top}\braceit{\tilde{P}_{t} - \tilde{P}_{t-1}}\fat{\zeta}
    \right]\nonumber\\
    &=\alpha^{\varrho}\sum_{t=2}^{T}\sum_{n=1}^{N}\text{Var}\left[
    [\tilde{P}_{t} - \tilde{P}_{t-1}]_{i}\fat{\zeta}
    \right]\nonumber\\
    &=\alpha^{\varrho}\sum_{t=2}^{T}\text{Tr}\left[
    \braceit{\tilde{P}_{t} - \tilde{P}_{t-1}}\widehat{\Sigma}\braceit{\tilde{P}_{t} - \tilde{P}_{t-1}}^{\top}
    \right],\label{eq:rho_var_ref}
\end{align}
\end{subequations}
which is a convex function. This derivation shows that term \eqref{eq:rho_var} minimizes variability, e.g., inter-temporal pressure ramps, and the variance of stochastic pressures simultaneously.  

\subsubsection{Network Topology Optimization}\label{subsubsec:topology_opt} Network topology has substantial impacts on disturbance propagation in natural gas networks. Fortunately, the topology can be optimized by the network operator by activating \textit{binary valves}. Unlike control valves with variable decompression rates $\tilde{\kappa}_{t}(\fat{\zeta}^{t})$, the binary valves decouple gas states at both sides of the valve \cite{fugenschuh2015chapter}. Hence, by activating binary valves, the network operator decouples network parts to minimize uncertainty propagation. To optimize topology, consider $V$ binary valves installed in the network and consider a set $\mathcal{C}=\{1,\dots,C\}$ of possible network topologies with $C=2^{V}$. To select a network topology, consider binary variables $v_{1},\dots,v_{C}\in\{0;1\}$: $v_{c}=1$ if topology $c$ is chosen, and  $v_{c}=0$ otherwise. Only one topology can be chosen for the entire control horizon, meaning 
\begin{subequations}\label{eq:topology_opt}
\begin{align}\label{eq:unique_topology}
  \sum_{c=1}^{C}v_{c} = 1.
\end{align}
Suppose that for each topology $c\in\mathcal{C}$, sensitivity coefficients $w_{0c},W_{1tc}$ and $W_{2tc}$ in \eqref{eq:wey_lin} and reference node pressure $\mathring{\varrho}_{\text{r}tc}$ are known to the network operator, i.e., they can always be pre-computed offline by solving a finite series of deterministic optimization problems. Then, stochastic gas flows proxy  $\tilde{\phi}_{tc}(\fat{\zeta}^{t})\in\mathbb{R}^{E}$ at time stage $t$ in network topology $c$ writes as
\begin{align}
    &\tilde{\phi}_{tc}(\fat{\zeta}^{t})=w_{0tc} + W_{1tc}\tilde{\varrho}_{t}(\fat{\zeta}^{t}) + W_{2tc}\tilde{\kappa}_{t}(\fat{\zeta}^{t}). 
\end{align}
Here, the sensitivities amount to zero for pipelines with binary valves when those are activated, hence preventing the gas flow. Then, the following constraints map the proxy gas flows and reference node pressure to the ones under the chosen topology: 
\begin{align}\label{eq:top_opt_last}
    &\tilde{\varphi}_{t}(\fat{\zeta}^{t})=\sum_{c=1}^{C}v_{c}\tilde{\phi}_{tc}(\fat{\zeta}^{t}),\;
    &\tilde{\varrho}_{\text{r}t}(\fat{\zeta}^{t})=\sum_{c=1}^{C}v_{c}\mathring{\varrho}_{\text{r}tc},
\end{align}
\end{subequations}
because only one topology can be selected as per equation \eqref{eq:unique_topology}. Then, stochastic gas flow variables $\tilde{\varphi}_{t}(\fat{\zeta}^{t})$ participate in the rest of the gas flow equations in \eqref{stoch_intr_as}. The tractable reformulation of stochastic equations \eqref{eq:topology_opt} is similar to that in Section \ref{subsubsec:ref_as}. Yet, the flow equation in \eqref{eq:top_opt_last} involves the product of binary and continuous variables, which is addressed by the standard Big-M approach. We provide the full problem reformulation as a mixed-integer SOCP (MISOCP) in \cite{dvorkin2021}. 

\section{Numerical Experiments}\label{sec:experiments}
\subsection{Experimental Setup}
We run experiments on the IEEE 118-bus system from \cite{lubin2019chance}, with 11 wind farms distributed among 3 zones, and on the 48-node natural gas network from \cite{dvorkin2020stochastic}. The two systems are coupled through 9 CGGT plants sited in the gas network as in Fig. \ref{1b}. The goal is to accommodate the disturbances from renewable generation in the gas network through a 5-stage control horizon with several hours between the stages. To model renewable uncertainty, we assume that the intertemporal changes of wind power generation in each zone is conditioned on the realization of the random variable $\boldsymbol{\zeta}$ with \textit{arbitrary} distribution but with known first- and second-order moments $\mu=\mathbb{1}$ and $\Sigma=\text{dg}[0,\sigma_{2}^{2},\dots,\sigma_{13}^{2}]$, respectively, with $\sigma_{i}^{2}=0.15,\forall i$. In terms of stochastic process \eqref{eq:renewable_process}, this means $k_{1}=1$, $k_{t}=3,\forall t\geqslant2$ and $k=13$. We then choose matrices $\Omega_{t}$ in \eqref{eq:renewable_process} to obtain the descending stochastic wind power generation process depicted in Fig. \ref{1a} (top). Then, following the procedure in Section \ref{subsec:from_wind_to_gas}, we solve a chance-constrained optimal power flow problem and convert the resulting stochastic CCGT generation into the ascending stochastic natural gas extraction process displayed in Fig. \ref{1a} (bottom). The standard deviation of stochastic gas extractions amounts to up to 7.2\% of their nominal values. When optimizing linear decision rules \eqref{eq:LDRs} to accommodate this stochastic process, we require the stochastic network state trajectories to remain within network limits with individual constraint satisfaction probability at least 99.5\% (i.e., $1-\overline{\varepsilon}_{t}=0.995$). Finally, we extract $10^3$ uncertainty samples from $
N(\mu,\Sigma)$ for the out-of-sample analysis.

Next, we provide selected results for control applications from Section \ref{subsec:control_applications}. Solution time does not exceed 1.2 sec for the SOCP and 24.5 sec for the MISOCP programs on average with the MOSEK solver on a standard laptop. All modeling data and codes to replicate the results are available in \cite{dvorkin2021}.  

\subsection{Results}

\begin{table*}[t]
\footnotesize
\caption{Expected cost, state variability and feasibility across the 5-stage control horizon under deterministic and stochastic control}
\label{tab:summary}
\centering
\setlength\tabcolsep{3.5pt}
\begin{tabular}{llcccccc}
\toprule
\multirow{3}{*}{Parameter} & \multirow{3}{*}{Unit} & \multirow{3}{*}{\begin{tabular}[c]{@{}c@{}}Deterministic \\ control policy\end{tabular}} & \multicolumn{5}{c}{Stochastic control policy}\\
\cmidrule(lr){4-8}
 &  &  & \multirow{2}{*}{\begin{tabular}[c]{@{}c@{}}Base\end{tabular}} & \multirow{2}{*}{\begin{tabular}[c]{@{}c@{}}Linepack-\\ agnostic\end{tabular}} &  \multicolumn{3}{c}{Variability-aware} \\
\cmidrule(lr){6-8}
 &  &  &  & & $\alpha^{\varrho}=10$ & $\alpha^{\varrho}=50$ & $\alpha^{\varrho}=100$\\
\midrule
Expected gas injection cost & \$1000 & 644.8 (\textbf{94.6\%}) & 681.7 (\textbf{100.0\%}) & 752.1 (\textbf{110.3\%}) & 694.4 (\textbf{101.9\%}) & 701.2 (\textbf{102.9\%}) & 703.5 (\textbf{103.2\%})\\
Pressure variability term $\frac{\eqref{eq:rho_var}}{\alpha^{\varrho}}$ & MPa & \textcolor{white}{0}72.01 (\textbf{189.7\%}) & \textcolor{white}{0}38.0 (\textbf{100.0\%}) & \textcolor{white}{0}66.0 (\textbf{173.6\%}) & \textcolor{white}{00}7.8 \textcolor{white}{0}(\textbf{20.5\%}) & \textcolor{white}{00}7.3 \textcolor{white}{0}(\textbf{19.2\%}) & \textcolor{white}{00}7.2 \textcolor{white}{0}(\textbf{19.1\%})\\
\midrule
\begin{tabular}{@{}l@{}}
Expected\;/\;worst-case magnitude\\
of pressure constraint violations
\end{tabular}
& MPa & 77.42\;/\;147.76 & 0.00\;/\;0.00 & 0.00\;/\;0.00 & 0.00\;/\;0.00 & 0.00\;/\;0.00 & 0.00\;/\;0.00\\
\cmidrule(lr){1-8}
\begin{tabular}{@{}l@{}}
Expected\;/\;worst-case magnitude\\ of gas mass constraint violations
\end{tabular} & MMSCFD & 26.86\;/\;146.91 & 0.01\;/\;0.02 & 0.13\;/\;0.13 & 0.01\;/\;0.02 & 0.02\;/\;0.02  & 0.02\;/\;0.02\\
\midrule
\begin{tabular}{@{}l@{}}First-stage gas injection\\ $\sum_{n\in\mathcal{N}}\vartheta_{1n}$\end{tabular} & MMSCFD & 2924.8 & 3229.0 & 3203.8 & 3233.2 & 3219.5 & 3219.1\\
\cmidrule(lr){1-8}
\begin{tabular}{@{}l@{}}Expected compressor deployment \\ $\sum_{\ell\in\mathcal{E}_{c}}^{t\in\mathcal{T}}\mathbb{E}_{\mathbb{P}_{\boldsymbol{\zeta}^{t}}}[\tilde{\kappa}_{t\ell}(\fat{\zeta}^{t})]$\end{tabular}  & kPa & 7127.9 & 10225.3  & 10912.7 & 12464.8 & 12451.1 & 12459.8\\
\cmidrule(lr){1-8}
\begin{tabular}{@{}l@{}}Expected valve deployment \\ $\sum_{\ell\in\mathcal{E}_{v}}^{t\in\mathcal{T}}\mathbb{E}_{\mathbb{P}_{\boldsymbol{\zeta}^{t}}}[\tilde{\kappa}_{t\ell}(\fat{\zeta}^{t})]$\end{tabular} & kPa & 0.0 & 714.8 & 1251.27 & 2281.0 & 2604.2 & 2647.6\\
\bottomrule
\end{tabular}
\vspace{-0.75cm}
\end{table*}

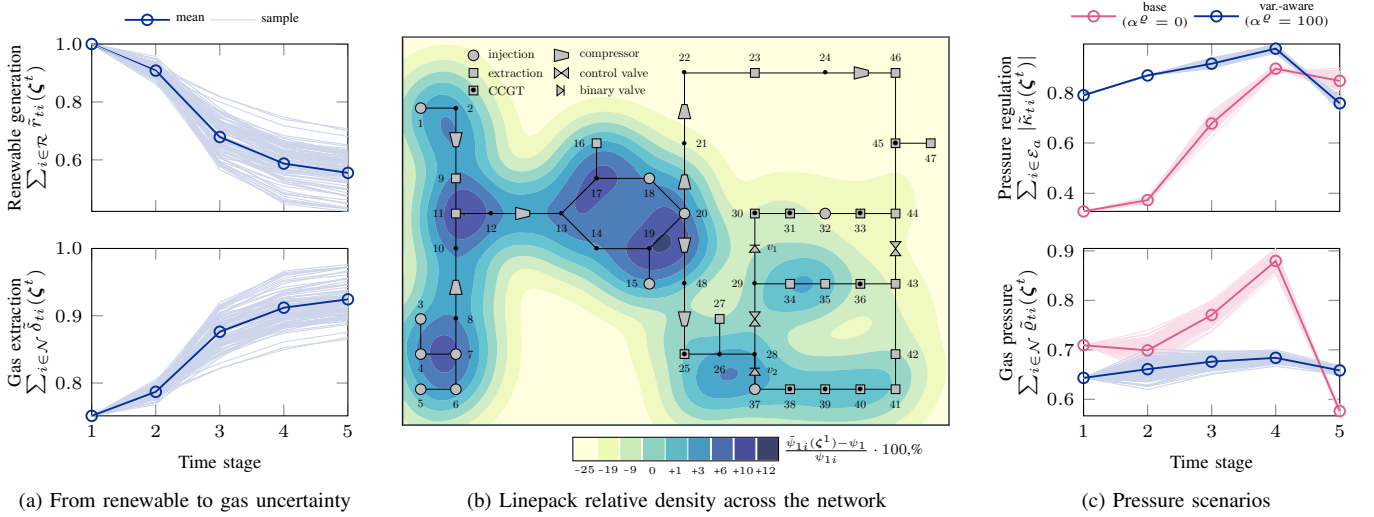
\begin{figure*}[t] 
    \centering
    \subfloat[From renewable to gas uncertainty\label{1a}]{%
      \begin{tikzpicture}[font=\scriptsize]
        \begin{axis}[width=0.275\linewidth,height=0.21\linewidth,
            ylabel = {Renewable generation \\ $\sum_{i\in\mathcal{R}}\tilde{r}_{ti}(\boldsymbol{\zeta}^{t})$},
            xtick={1,2,3,4,5},
            xticklabels={,,},
            y tick label style={
                        /pgf/number format/.cd,
                        fixed,
                        fixed zerofill,
                        precision=1,
                        /tikz/.cd
            },
            y label style={align=center,yshift=-0.15cm},
            enlargelimits=false,
            legend style={cells={align=center},legend columns=-1,font=\tiny,draw=none,xshift=-0.4cm,yshift=0.6cm}
            ]
            \foreach \a in {2,3,...,100}{
            \addplot [color=blue!20,mark=none, forget plot]table [x index=0, y index=\a, col sep=comma] {PSCC_2022_code/results/wind_results.csv};
            }
            \addplot[color=blue,mark=o,line width = 0.025cm] table [x index=0, y index=1, col sep=comma] {PSCC_2022_code/results/wind_results.csv};
            \addlegendentry{mean};
            \addplot [color=blue!20,mark=none]table [x index=0, y index=3, col sep=comma] {PSCC_2022_code/results/wind_results.csv};
            \addlegendentry{sample};
        \end{axis}
        
        \begin{axis}[width=0.275\linewidth,height=0.21\linewidth,yshift = -0.15\linewidth,
                    ylabel = {Gas extraction\\ $\sum_{i\in\mathcal{N}}\tilde{\delta}_{ti}(\boldsymbol{\zeta}^{t})$},
                    xtick={1,2,3,4,5},
                    ymax = 1.0,
                    xlabel = {Time stage},
                    y tick label style={
                        /pgf/number format/.cd,
                        fixed,
                        fixed zerofill,
                        precision=1,
                        /tikz/.cd
                    },
                    enlargelimits=false,
                    y label style={align=center,yshift=-0.15cm}
                    ]
            \foreach \a in {2,3,...,100}{
            \addplot [color=blue!20,mark=none]table [x index=0, y index=\a, col sep=comma] {PSCC_2022_code/results/gas_coms_results.csv};
            }
            \addplot[color=blue,mark=o,line width = 0.025cm] table [x index=0, y index=1, col sep=comma] {PSCC_2022_code/results/gas_coms_results.csv};
        \end{axis}
      \end{tikzpicture}
      }
    \hfill
    \subfloat[Linepack relative density across the network \label{1b}]{%
    
        \begin{tikzpicture}
        \node[inner sep=0pt] (a) at (0,0) {\includegraphics[width=0.4\linewidth, frame]{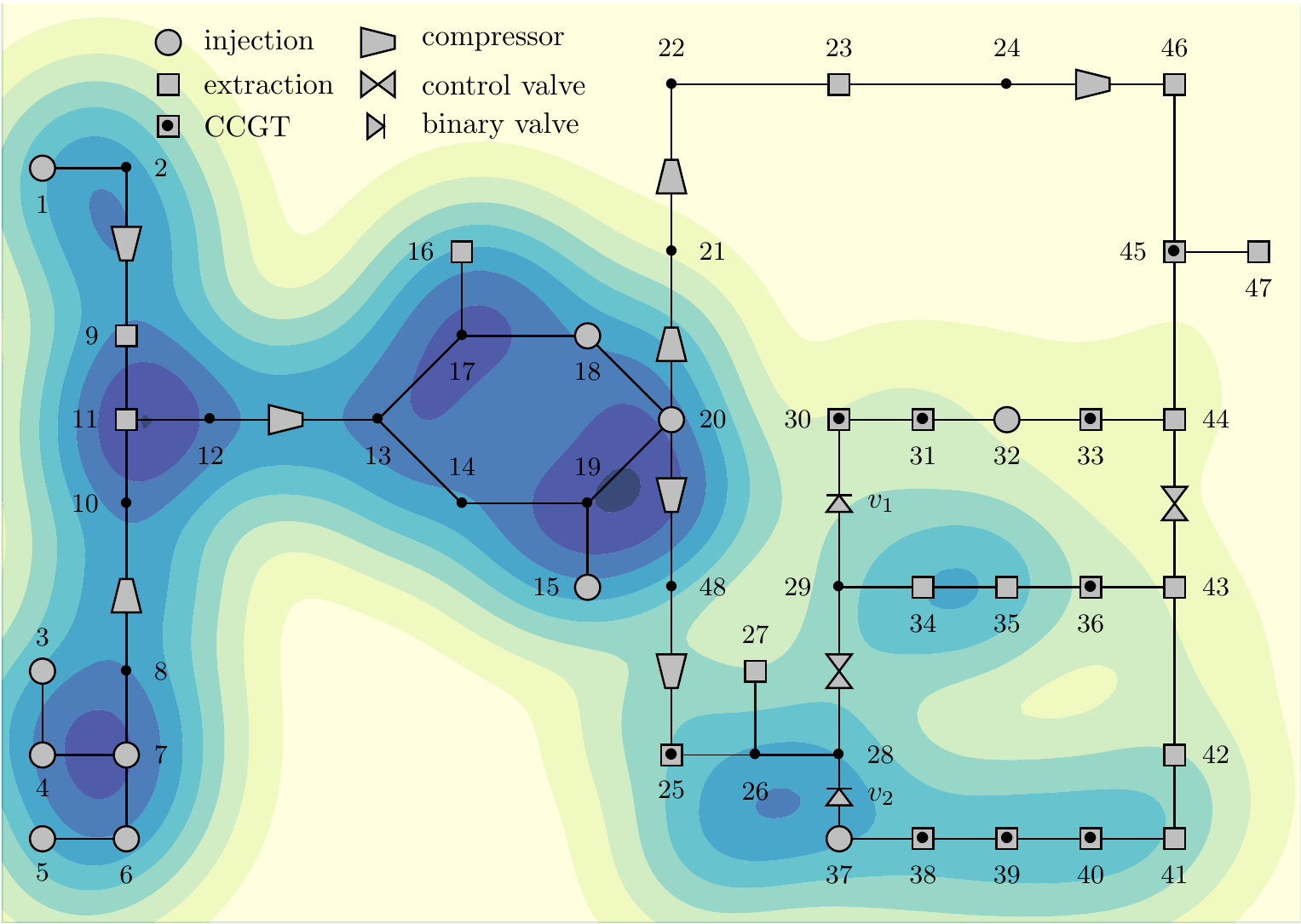}};
        \begin{scope}
        \node [below of = a, yshift = -1.9cm] (l) {\includegraphics[width=0.15\linewidth, frame]{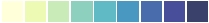}};
        \node [above of = l, yshift = -1.28cm,font=\scriptsize, scale = 0.8, xshift = -1.515cm,font=\tiny] {$\shortminus25$};
        \node [above of = l, yshift = -1.28cm,font=\scriptsize, scale = 0.8, xshift = -1.125cm,font=\tiny] {$\shortminus19$};
        \node [above of = l, yshift = -1.28cm,font=\scriptsize, scale = 0.8, xshift = -0.76cm,font=\tiny] {$\shortminus9$};
        \node [above of = l, yshift = -1.28cm,font=\scriptsize, scale = 0.8, xshift = -0.385cm,font=\tiny] {$0$};
        \node [above of = l, yshift = -1.28cm,font=\scriptsize, scale = 0.8, xshift = -0.0cm,font=\tiny] {+$1$};
        \node [above of = l, yshift = -1.28cm,font=\scriptsize, scale = 0.8, xshift = 0.375cm,font=\tiny] {+$3$};
        \node [above of = l, yshift = -1.28cm,font=\scriptsize, scale = 0.8, xshift = 0.750cm,font=\tiny] {+$6$};
        \node [above of = l, yshift = -1.28cm,font=\scriptsize, scale = 0.8, xshift = 1.125cm,font=\tiny] {+$10$};
        \node [above of = l, yshift = -1.28cm,font=\scriptsize, scale = 0.8, xshift = 1.500cm,font=\tiny] {+$12$};
        \node [above of = l, yshift = -1.0cm,font=\scriptsize, scale = 0.8, xshift = 2.95cm,font=\scriptsize] {$\frac{\tilde{\psi}_{1i}(\boldsymbol{\zeta}^{1}) - \psi_{1}}{\psi_{1i}}\cdot100$,\%};
        \end{scope}
        \end{tikzpicture}
        }
    \hfill   
    \subfloat[Pressure scenarios\label{1c}]{%
        \begin{tikzpicture}[font=\scriptsize]
        \begin{axis}[width=0.275\linewidth,height=0.21\linewidth,
                    ylabel = {Pressure regulation \\ $\sum_{i\in\mathcal{E}_{a}}|\tilde{\kappa}_{ti}(\boldsymbol{\zeta}^{t})|$},
                    xtick={1,2,3,4,5},
                    xticklabels={,,},
                    y tick label style={
                        /pgf/number format/.cd,
                        fixed,
                        fixed zerofill,
                        precision=1,
                        /tikz/.cd
                    },
                    legend style={cells={align=center,fill=none},legend columns=-1,font=\tiny,draw=none,xshift=0.1cm,yshift=0.75cm},
                    enlargelimits=false,
                    y label style={align=center,yshift=-0.15cm}
                    ]
            \foreach \a in {2,3,...,100}{
            \addplot [color=red!20,mark=none, forget plot]table [x index=0, y index=\a, col sep=comma] {PSCC_2022_code/results/sto_pres_reg_wo_var.csv};
            }
            \foreach \a in {2,3,...,100}{
            \addplot [color=blue!20,mark=none, forget plot]table [x index=0, y index=\a, col sep=comma] {PSCC_2022_code/results/sto_pres_reg_wt_var.csv};
            }
            \addplot[color=red,mark=o,line width = 0.025cm] table [x index=0, y index=1, col sep=comma] {PSCC_2022_code/results/sto_pres_reg_wo_var.csv};
            \addlegendentry{base \\ ($\alpha^{\varrho}=0$)};
            \addplot[color=blue,mark=o,line width = 0.025cm] table [x index=0, y index=1, col sep=comma] {PSCC_2022_code/results/sto_pres_reg_wt_var.csv};
            \addlegendentry{var.-aware \\ ($\alpha^{\varrho}=100$)};
        \end{axis}
        \begin{axis}[width=0.275\linewidth,height=0.21\linewidth,yshift = -0.15\linewidth,
                    ylabel = {Gas pressure\\ $\sum_{i\in\mathcal{N}}\tilde{\varrho}_{ti}(\boldsymbol{\zeta}^{t})$},
                    xtick={1,2,3,4,5},
                    xlabel = {Time stage},
                    y tick label style={
                        /pgf/number format/.cd,
                        fixed,
                        fixed zerofill,
                        precision=1,
                        /tikz/.cd
                    },
                    enlargelimits=false,
                    y label style={align=center,yshift=-0.15cm}
                    ]
            \foreach \a in {2,3,...,100}{
            \addplot [color=red!20,mark=none]table [x index=0, y index=\a, col sep=comma] {PSCC_2022_code/results/sto_pressure_wo_var.csv};
            }
            \foreach \a in {2,3,...,100}{
            \addplot [color=blue!20,mark=none]table [x index=0, y index=\a, col sep=comma] {PSCC_2022_code/results/sto_pressure_wt_var.csv};
            }
            \addplot[color=red,mark=o,line width = 0.025cm] table [x index=0, y index=1, col sep=comma] {PSCC_2022_code/results/sto_pressure_wo_var.csv};
            \addplot[color=blue,mark=o,line width = 0.025cm] table [x index=0, y index=1, col sep=comma] {PSCC_2022_code/results/sto_pressure_wt_var.csv};
        \end{axis}
      \end{tikzpicture}
        }
  \caption{Multi-stage natural gas network optimization: (a) Normalized total renewable power generation scenarios (top) and their conversion into normalized total gas extraction (bottom) as explained in Section \ref{subsec:from_wind_to_gas}, (b) Density plot of the relative difference between the deterministic ($\psi_{1i}$) and the base stochastic ($\tilde{\psi}_{1i}
  (\boldsymbol{\zeta}^{1})$) first-stage linepack decisions for $i=1,\dots,E$, (c) Normalized scenarios of total pressure regulation (top) and aggregated nodal pressures (bottom) under the base stochastic (variability-agnostic) and stochastic variability-aware network control policies.}
  \label{fig1} 
  \vspace{-0.5cm}
\end{figure*}

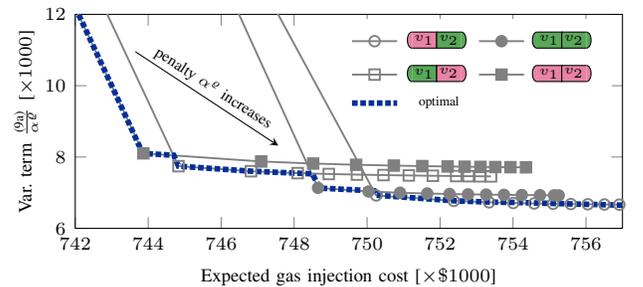
\begin{figure}
    \centering
        \begin{tikzpicture}[font=\scriptsize]
        \begin{axis}[width=1\linewidth,height=0.5\linewidth,
            xlabel = {Expected gas injection cost [$\times\$1000$]},
            ylabel = {Var. term $\frac{\eqref{eq:rho_var}}{\alpha^{\varrho}}$ [$\times1000$]},
            ymin = 6,
            ymax = 12,
            xmin = 742,
            xmax = 757,
            y label style={align=center,yshift=-0.15cm},
            legend style={scale=0.8,cells={align=center},legend columns=2,font=\tiny,draw=none,xshift=0.cm,yshift=0.cm}
            ]
            \addplot[color=gray,mark=o,line width = 0.025cm] table [x=cost, y=std, col sep=comma] {PSCC_2022_code/results/var_trade_off_conf_c_1.csv};
            \addlegendentry{
            \begin{tikzpicture}
            \draw[fill=red!75] (0,0) edge[bend left=90,fill=red!75] (0,0.075) -- (0.045,0.0)  -- (0.045,0.075) -- (0.0,0.075) node[pos=0.5,yshift=-2.6] {$v_{1}$};
            \draw[fill=green!75] (0.09,0.075) edge[bend left=90,fill=green!75] (0.09,0.0) -- (0.045,0.075) -- (0.045,0.0) -- (0.09,0.0) node[pos=0.5,yshift=3.5] {$v_{2}$};
            \end{tikzpicture}
            }
            \addplot[color=gray,mark=*,line width = 0.025cm] table [x=cost, y=std, col sep=comma] {PSCC_2022_code/results/var_trade_off_conf_c_2.csv};
            \addlegendentry{
            \begin{tikzpicture}
            \draw[fill=green!75] (0,0) edge[bend left=90,fill=green!75] (0,0.075) -- (0.045,0.0)  -- (0.045,0.075) -- (0.0,0.075) node[pos=0.5,yshift=-2.6] {$v_{1}$};
            \draw[fill=green!75] (0.09,0.075) edge[bend left=90,fill=green!75] (0.09,0.0) -- (0.045,0.075) -- (0.045,0.0) -- (0.09,0.0) node[pos=0.5,yshift=3.5] {$v_{2}$};
            \end{tikzpicture}
            }
            \addplot[color=gray,mark=square,line width = 0.025cm] table [x=cost, y=std, col sep=comma] {PSCC_2022_code/results/var_trade_off_conf_c_3.csv};
            \addlegendentry{
            \begin{tikzpicture}
            \draw[fill=green!75] (0,0) edge[bend left=90,fill=green!75] (0,0.075) -- (0.045,0.0)  -- (0.045,0.075) -- (0.0,0.075) node[pos=0.5,yshift=-2.6] {$v_{1}$};
            \draw[fill=red!75] (0.09,0.075) edge[bend left=90,fill=red!75] (0.09,0.0) -- (0.045,0.075) -- (0.045,0.0) -- (0.09,0.0) node[pos=0.5,yshift=3.5] {$v_{2}$};
            \end{tikzpicture}
            }
            \addplot[color=gray,mark=square*,line width = 0.025cm] table [x=cost, y=std, col sep=comma] {PSCC_2022_code/results/var_trade_off_conf_c_4.csv};
            \addlegendentry{
            \begin{tikzpicture}
            \draw[fill=red!75] (0,0) edge[bend left=90,fill=red!75] (0,0.075) -- (0.045,0.0)  -- (0.045,0.075) -- (0.0,0.075) node[pos=0.5,yshift=-2.6] {$v_{1}$};
            \draw[fill=red!75] (0.09,0.075) edge[bend left=90,fill=red!75] (0.09,0.0) -- (0.045,0.075) -- (0.045,0.0) -- (0.09,0.0) node[pos=0.5,yshift=3.5] {$v_{2}$};
            \end{tikzpicture}
            }
            \addplot[color=blue,densely dotted,mark=none,line width = 0.075cm] coordinates{(742,12.2)(743,10)(743.85,8.1)(744.7,8.05)(744.83,7.74375)(746.817,7.60241)(748.5,7.5287)(748.7,7.1287)(750.2,7.051)(750.25,6.93)(752.442,6.77)(753.4,6.735)(754,6.728)(755,6.7)(756.99,6.65)};
            \addlegendentry{optimal}
        \end{axis}
        \draw[,->,>=stealth] (0.8,2.35) -- node[rotate=-33,above,pos=0.5,scale=0.85,yshift=-0.05cm] {penalty $\alpha^{\varrho}$ increases} (2.7,1.1);
        \end{tikzpicture}
        \caption{Trade-offs between the expected cost and pressure variability under four network topologies for varying penalty factor $\alpha^{\varrho}$. The green/red color indicates that the binary valve is activated/deactivated. The optimal frontier (depicted in dashed blue) is obtained by solving the mixed-integer topology optimization problem from Section \ref{subsubsec:topology_opt} for various assignments $\alpha^{\varrho}$. The variance of renewable generation is increased in this experiment from 0.15 to 0.20 for more illustrative results.}
        \label{fig:topology_opt}
\end{figure}

We first compare deterministic and stochastic control policies in Table \ref{tab:summary}. The deterministic policy is obtained by solving problem \eqref{prob:stoch_tr} when replacing chance constraints \eqref{prob_cc_soc_min_flow}--\eqref{prob_cc_soc_last} with their deterministic counterparts to constraint the nominal variable components only. The base stochastic policy is optimized following the application in Section \ref{subsubsec:linepack}: we set factor $\alpha^{\vartheta}=2.5\%$ to minimize the adjustment of gas injections, while exploiting the linepack as the main source of flexibility. The results in the $3^{\text{rd}}$ and $4^{\text{th}}$ columns of Table \ref{tab:summary} show that the expected cost under deterministic policy is 5.4\% smaller than under the base stochastic policy. However, the out-of-sample analysis of the deterministic solution on 1000 uncertainty realization scenarios demonstrates substantial magnitudes of pressure and gas mass constraint violation across the 5-stage horizon, both in expectation and in the average worst-case scenario (computed below the $5\%$-quantile of the violation magnitude distribution). Improving on the deterministic solution, the base stochastic control policy demonstrates practically zero constraint violation at the expense of a 5.4\% increase in expected cost. The difference between the two policies is that the stochastic one requires more nominal gas injection at the first stage, which is then stored as a linepack in the western network part to be deployed later to balance uncertain extractions, as shown by the density plot in Fig. \ref{1b}. 

Moreover, optimizing the base stochastic policy with the exact reformulation of the double-sided chance constraints in \eqref{prob:stoch_tr} demonstrates a less conservative solution than the Chebyshev approximation (as in \cite{ratha2020affine} and \cite{dvorkin2020stochastic}), which guarantees the same constraint satisfaction probability but at larger expected cost of \$698,643, which is 2.5\% more expensive than in Table \ref{tab:summary}. 

Next, we estimate the cost-saving potential of linepack flexibility in balancing stochastic gas extractions. To do so, we compare the base policy with the linepack-agnostic policy which disregards the flexibility of linepack. To obtain this policy, we reduce the recourse (adjustment) of linepack stochastic variable $\tilde{\psi}_{t}(\fat{\zeta}^{t})$ by introducing the following constraint:
\begin{align}
    &\norm{\widehat{F}[\Psi_{t}]_{\ell}S_{t}}
    \leqslant \alpha^{\psi}[\Psi_{t}]_{\ell}S_{t}\widehat{\mu},\; \forall \ell\in\mathcal{E}, \forall t\in\mathcal{T},
\end{align}
to limit the standard deviation of linepack variable by an $\alpha^{\psi}-$portion of its nominal value. We empirically find out that the minimal value of $\alpha^{\psi}$ that ensures the feasibility of chance-constrained program \eqref{prob:stoch_tr} is 4\%, for which we report the results in the $5^{\text{th}}$ column of Table \ref{tab:summary}. Observe, that the feasibility performance of the base and linepack-agnostic strategies remain similar due to guarantees of chance constraints, but the ignorance to linepack flexibility results in a substantial increase in the expected cost. Here, we estimate the cost-saving potential of the linepack flexibility at 10.3\%.

Next, consider three variability-aware policies for varying pressure variability penalty factor $\alpha^{\varrho}$ in Table \ref{tab:summary}. Even with a small penalty, the pressure variability already drops to 20.5\% relative to the base policy at a small increase in operating costs by 1.9\%. An increasing penalty further reduces the variability but at a larger cost, thus enabling the network operator to trade off between the expected cost and network state variability. Table \ref{tab:summary} and Figure \ref{1c} (top) illustrate that variability-aware policies leverage active pipelines more intensively, with more notable use of control valves to regulate pressures. This yields substantially less variable and uncertain nodal pressure profiles, as shown in Fig. \ref{1c} (bottom). 

Finally, we analyze the  network topology optimization from Section \ref{subsubsec:topology_opt}. The topology is optimized by switching two binary valves $v_{1}$ and $v_{2}$ in edges (29,30) and (28,37), respectively, to provide the optimal cost-variability frontier to the network operator, as depicted in Fig. \ref{fig:topology_opt}. For example, the expected operating cost of $\approx\$754\times10^3$ is achieved under the original topology with no activated valve, as well as under two other topologies with the second or both valves active. However, it is variability-optimal to maintain this cost level by activating the second valve to reduce the pressure variability measure by 13.2\%, from 7.71 to 6.69 MPa.

\section{Conclusions}\label{sec:conclusion}

To address propagating uncertainty and variability of renewable generation, we developed multi-stage stochastic control policies for natural gas network operators to leverage the linepack flexibility. Using chance-constrained programming, these policies were optimized for compressors and valves to satisfy feasibility and variability criteria without substantial deviations from the nominal natural gas supply injections. 

Using a realistically sized natural gas network, we estimated the cost-saving potential of linepack flexibility at 10.3\% of operating costs in expectation. We also found out that the system state variability can be reduced to 20.5\% of the nominal value at a small 1.9\% increase in expected costs. Finally, we concluded that the network topology must be optimized itself to enable new cost-variability trade offs to network operators. 

Future research will focus on how the proposed stochastic control policies agree with a dynamic natural gas modeling to explicitly account for the time delays in gas transportation and linepack storage across the multi-stage control horizon.

\bibliographystyle{IEEEtran}
\bibliography{references.bib}
\endgroup


\section*{Supplementary material}

# Multi-Stage Linear Decision Rules for Stochastic Control of Natural Gas Networks with Linepack (Online Appendix)

Vladimir Dvorkin*, Dharik Mallapragada*, Audun Botterud*, Jalal Kazempour† and Pierre Pinson†
*Massachusetts Institute of Technology, Cambridge, MA, USA.
†Technical University of Denmark, Kng. Lyngby, Denmark.

## I. Multi-Stage Optimal Power Flow Problem Formulation

We model power system as a network with $N$ nodes and $\Lambda$ transmission lines connecting those nodes. At any stage $t$ of the control horizon, power generation from conventional and renewable units must satisfy $L$ loads collected in vector $d_t \in \mathbb{R}_+^L$ and allocated in the network according to incidence matrix $M_d \in \mathbb{R}^{N \times L}$. Towards this goal, the power system operator dispatches $P$ conventional and $R$ renewable generation units, which are allocated in the network according to incidence matrices $M_g \in \mathbb{R}^{N \times P}$ and $M_r \in \mathbb{R}^{N \times R}$, respectively, while meeting the minimum/maximum generation limits $\underline{g}/\overline{g} \in \mathbb{R}_+^P$ of conventional generators and power flow limits $\overline{f} \in \mathbb{R}^\Lambda$. The cost of conventional generation is modeled using a linear function with coefficients collected in vector $c \in \mathbb{R}_+^P$. The power flows are modelled using the *DC* power flow approximation and the matrix of power transfer distribution factors $F \in \mathbb{R}^{N \times \Lambda}$. Given the stochastic renewable in-feed $\tilde{r}_t(\zeta^t)$ from Section II-C, the multi-stage optimal power flow problem is addressed as the following stochastic program:

$$\underset{\tilde{g}_t}{\text{minimize}} \quad \mathbb{E}_{\mathbb{P}_\zeta}\left[\sum_{t=1}^T c^\top \tilde{g}_t(\zeta^t)\right] \tag{1a}$$

$$\text{subject to:} \quad \mathbb{P}_{\zeta^t}\left[\mathbb{1}^\top (M_g \tilde{g}_t(\zeta^t) + M_r \tilde{r}_t(\zeta^t) - M_d d_t) = 0\right] = 1, \tag{1b}$$

$$\mathbb{P}_{\zeta^t}\begin{bmatrix} -\overline{f} \leqslant F(M_g \tilde{g}_t(\zeta^t) + M_r \tilde{r}_t(\zeta^t) - M_d d_t) \leqslant \overline{f}, \\ \underline{g} \leqslant \tilde{g}_t(\zeta^t) \leqslant \overline{g}, \end{bmatrix} \geqslant 1 - \varepsilon_t, \quad \forall t \in \mathcal{T}, \tag{1c}$$

where the objective function to be minimized is the expected generation costs, subject to two sets of probabilistic constraints. Constraint (1b) is enforced to preserve the power balance at each time stage with probability 1, and constraint (1c) is enforced to satisfy the power flow and generation limits jointly at every time stage with probability at least $1 - \varepsilon_t$, where $\varepsilon_t$ is a small parameter chosen by the power system operator. Defining generator response as the following multi-stage linear decision rule:

$$\tilde{g}_t(\zeta^t) = G_t S_t \zeta, \quad \forall t \in \mathcal{T}, \tag{2}$$

with matrix $G \in \mathbb{R}^{N \times k^t}$ of variables, and taking the same reformulation steps as for the gas network optimization problem in Section III-A, this chance-constrained program reformulates into a computationally tractable SOCP program:

$$\underset{G_t, y_t^f, x_t^f}{\text{minimize}} \quad \sum_{t=1}^T c^\top G_t S_t \widehat{\mu} \tag{3a}$$

$$\text{subject to:} \quad \mathbb{1}^\top (M_g G_t + M_r \Omega_t - M_d \begin{bmatrix} d_t & \mathbb{0}_{L \times k^t - 1} \end{bmatrix}) = \mathbb{0}, \tag{3b}$$

$$\sqrt[-2]{\overline{\varepsilon}_t} \left\| \begin{matrix} \widehat{F}[(M_g G_t + M_r \Omega_t - M_d \begin{bmatrix} d_t & \mathbb{0}_{L \times k^t - 1} \end{bmatrix}) S_t]_\ell^\top \\ y_{t\ell}^f \end{matrix} \right\| \leqslant \overline{f}_\ell - x_{t\ell}^f, \tag{3c}$$

$$\sqrt[-2]{\overline{\varepsilon}_t} \left\| \begin{matrix} \widehat{F}[G_t S_t]_i^\top \\ y_{ti}^g \end{matrix} \right\| \leqslant \tfrac{1}{2}(\overline{g}_i - \underline{g}_i) - x_{ti}^g, \tag{3d}$$

$$\left|[(M_g G_t + M_r \Omega_t - M_d \begin{bmatrix} d_t & \mathbb{0}_{L \times k^t - 1} \end{bmatrix}) S_t]_\ell \widehat{\mu}\right| \leqslant y_{t\ell}^f + x_{t\ell}^f, \tag{3e}$$

$$\left|[G_t S_t]_i \widehat{\mu} - \tfrac{1}{2}(\overline{g}_i - \underline{g}_i)\right| \leqslant y_{ti}^g + x_{ti}^g, \tag{3f}$$

$$\overline{f}_\ell \geqslant x_{t\ell}^f \geqslant 0, \; y_{t\ell}^f \geqslant 0, \tag{3g}$$

$$\tfrac{1}{2}(\overline{g}_i - \underline{g}_i) \geqslant x_{ti}^g \geqslant 0, \; y_{ti}^g \geqslant 0, \tag{3h}$$

$$\forall t = 1, \ldots, T, \; \forall i = 1, \ldots, N, \; \forall \ell = 1, \ldots, \Lambda, \tag{3i}$$


where $\widehat{\mu}$ is the mean vector of renewable power deviations and $\widehat{F}$ is the Cholesky decomposition of covariance matrix $\widehat{\Sigma}$ of renewable power deviations. Here, linear constraint (3b) is the reformulation of (1b), and the set of conic and linear constraints (3c)–(3h) is a distributionally robust reformulation of the joint chance constraint (1c) on power flows and generation limits, where $\overline{\varepsilon}_t$ is the individual constraint violation probability. The power system operator optimizes variable $G$ of generator response (2) using program (3) and then converts this information into the stochastic gas extraction as explained in Section II-C.

## II. Network Topology Optimization Problem

To obtain a MISOCP reformulation of the network topology optimization problem from Section III-B3, conider the following steps. First, the linearized gas flow equation

$$\tilde{\phi}_{tc}(\zeta^t) = w_{0tc} + W_{1tc}\tilde{\varrho}_t(\zeta^t) + W_{2tc}\tilde{\kappa}_t(\zeta^t), \quad \forall t \in \mathcal{T}, \forall c \in \mathcal{C}, \tag{4a}$$

is formulated in multi-stage linear decision rules as

$$F_{tc}S_t\zeta = w_{0tc} + W_{1tc}P_tS_t\zeta + W_{2tc}K_tS_t\zeta, \quad \forall t \in \mathcal{T}, \forall c \in \mathcal{C}, \tag{4b}$$

where $F_{tc} \in \mathbb{R}^{E \times k^t}$ is the finite matrix of flow variable coefficients subject to optimization. Then, taking the path outlined in Section III-A2, this reformulates into the following set of linear constraints

$$\left(F_{tc} - [w_{0tc} \quad \mathbb{0}_{E \times (k^t-1)}] - W_{1tc}P_t - W_{2tc}K_t\right)S_t = \mathbb{0}, \quad \forall t \in \mathcal{T}, \forall c \in \mathcal{C}. \tag{4c}$$

The stochastic reference node pressure equality in (10c) reformulates similarly as before:

$$\tilde{\varrho}_{rt}(\zeta^t) = \sum_{c=1}^{C} v_c \mathring{\varrho}_{rtc} \iff \left([P_t]_r - \sum_{c=1}^{C} v_c [\mathring{\varrho}_{rtc} \quad \mathbb{0}_{1 \times (k^t-1)}]\right)S_t = \mathbb{0}, \quad \forall t \in \mathcal{T}. \tag{5}$$

Then, the stochastic gas flow constraint in (10c) is first reformulated as a set of linear equations:

$$\tilde{\varphi}_t(\zeta^t) = \sum_{c=1}^{C} v_c \tilde{\phi}_{tc}(\zeta^t)$$

$$\iff \Phi_t S_t \zeta = \sum_{c=1}^{C} v_c F_{tc} S_t \zeta$$

$$\stackrel{Z_{tc}=v_c F_{tc}}{\iff} \Phi_t S_t \zeta = \sum_{c=1}^{C} Z_{tc} S_t \zeta$$

$$\iff (\Phi_t - \sum_{c=1}^{C} Z_{tc})S_t = \mathbb{0}, \quad \forall t \in \mathcal{T}, \tag{6a}$$

where variable $Z_{tc} \in \mathbb{R}^{E \times k^t}$ substitutes the product of binary variable $v_c$ and continuous variable $F_{tc}$. This bilinear term is next reformulated using the standard Big-M approach:

$$Z_{tc}^{ij} = v_c F_{tc}^{ij} \iff \begin{cases} \underline{M} \leqslant Z_{tc}^{ij} \leqslant \overline{M} \\ \underline{M}v_c \leqslant Z_{tc}^{ij} \leqslant \overline{M}v_c \\ F_{tc}^{ij} - (1-v_c)\overline{M} \leqslant Z_{tc}^{ij} \leqslant F_{tc}^{ij} - (1-v_c)\underline{M} \end{cases} \quad \forall t \in \mathcal{T}, \forall c \in \mathcal{C}, \forall i = 1, \ldots, E, \forall j = 1, \ldots, k^t, \tag{6b}$$

where $Z_{tc}^{ij}$ (and $F_{tc}^{ij}$) denotes the entry of matrix $Z_{tc}$ (and $F_{tc}$) at position $(i,j)$, and $\overline{M}$ and $\underline{M}$ are sufficiently large positive and negative constants, respectively. Then, we obtain the following tractable MISOCP reformulation:

$$\text{minimize}_{\Theta_t, P_t, K_t, \Psi_t, \Phi_t^{(\cdot)}, Z_{tc}, F_{tc}, x_t^{(\cdot)}, y_t^{(\cdot)}, v_c}$$

$$\sum_{t=1}^{T}\left(c_1^\top \Theta_t S_t \widehat{\mu} + \text{Tr}\left[\Theta_t^\top \text{dg}[c_2]\Theta_t S_t \left(\widehat{\Sigma}+\widehat{\mu}\widehat{\mu}^\top\right)S_t^\top\right] + \alpha^\varrho \sum_{t=2}^{T} \text{Tr}\left[\left(\tilde{P}_t - \tilde{P}_{t-1}\right)\widehat{\Sigma}\left(\tilde{P}_t - \tilde{P}_{t-1}\right)^\top\right]\right) \tag{7a}$$

subject to:

$$(A^+\Phi_t^+ + A^-\Phi_t^- - \Theta_t S_t + BK_t + \Delta_t)S_t = \mathbb{0} \tag{7b}$$

$$\left(F_{tc} - [w_{0tc} \quad \mathbb{0}_{E \times (k^t-1)}] - W_{1tc}P_t - W_{2tc}K_t\right)S_t = \mathbb{0}, \quad \forall c \in \mathcal{C} \tag{7c}$$

$$\left([P_t]_r - \sum_{c=1}^{C} v_c [\mathring{\varrho}_{rtc} \quad \mathbb{0}_{1 \times (k^t-1)}]\right)S_t = \mathbb{0} \tag{7d}$$

$$(\Phi_t - \sum_{c=1}^{C} Z_{tc})S_t = \mathbb{0} \tag{7e}$$

$$\left(\Phi_t - \tfrac{1}{2}\Phi_t^+ - \tfrac{1}{2}\Phi_t^-\right) S_t = \mathbb{0} \tag{7f}$$

$$\begin{cases} \underline{M} \leqslant Z_{tc}^{ij} \leqslant \overline{M} \\ \underline{M}v_c \leqslant Z_{tc}^{ij} \leqslant \overline{M}v_c \\ F_{tc}^{ij} - (1-v_c)\overline{M} \leqslant Z_{tc}^{ij} \leqslant F_{tc}^{ij} - (1-v_c)\underline{M} \end{cases} \forall c \in \mathcal{C}, \forall i = 1,\ldots,E, \forall j = 1,\ldots,k^t \tag{7g}$$

$$\sum_{c=1}^{C} v_c = 1, \quad v_c \in \{0,1\} \tag{7h}$$

$$\left(\Psi_t - \tfrac{1}{2}\mathrm{dg}[s]\left(K_t + |A|^\top P_t\right)\right) S_t = \mathbb{0} \tag{7i}$$

$$(\Psi_t S_t - \Phi_t^+ + \Phi_t^-) S_t - \Psi_{(t\text{-}1)} S_{(t\text{-}1)} = \mathbb{0} \tag{7j}$$

$$\sqrt{\tfrac{1-\overline{\varepsilon}_t}{\overline{\varepsilon}_t}} \left\|\widehat{F}[\Phi_t S_t]_\ell^\top\right\| \leqslant [\Phi_t S_t \widehat{\mu}]_\ell \,^* \tag{7k}$$

$$\sqrt{\tfrac{1-\overline{\varepsilon}_t}{\overline{\varepsilon}_t}} \left\|\widehat{F}[\Psi_t S_t]_\ell^\top\right\| \leqslant [\Psi_t S_t \widehat{\mu} - \psi_0]_\ell \tag{7l}$$

$$\sqrt[-2]{\overline{\varepsilon}} \left\|\begin{matrix}\widehat{F}[\Theta_t S_t]_n^\top \\ y_{tn}^\vartheta\end{matrix}\right\| \leqslant \tfrac{1}{2}\left(\overline{\vartheta}_n - \underline{\vartheta}_n\right) - x_{tn}^\vartheta \tag{7m}$$

$$\sqrt[-2]{\overline{\varepsilon}} \left\|\begin{matrix}\widehat{F}[K_t S_t]_\ell^\top \\ y_{t\ell}^\kappa\end{matrix}\right\| \leqslant \tfrac{1}{2}\left(\overline{\kappa}_\ell - \underline{\kappa}_\ell\right) - x_{t\ell}^\kappa \tag{7n}$$

$$\sqrt[-2]{\overline{\varepsilon}} \left\|\begin{matrix}\widehat{F}[P_t S_t]_n^\top \\ y_{tn}^\varrho\end{matrix}\right\| \leqslant \tfrac{1}{2}\left(\overline{\varrho}_n - \underline{\varrho}_n\right) - x_{tn}^\varrho \tag{7o}$$

$$\left|[\Theta_t S_t]_n \widehat{\mu} - \tfrac{1}{2}\left(\overline{\vartheta}_n - \underline{\vartheta}_n\right)\right| \leqslant y_{tn}^\vartheta + x_{tn}^\vartheta \tag{7p}$$

$$\left|[K_t S_t]_\ell \widehat{\mu} - \tfrac{1}{2}\left(\overline{\kappa}_\ell - \underline{\kappa}_\ell\right)\right| \leqslant y_{t\ell}^\kappa + x_{t\ell}^\kappa \tag{7q}$$

$$\left|[P_t S_t]_n \widehat{\mu} - \tfrac{1}{2}\left(\overline{\varrho}_n - \underline{\varrho}_n\right)\right| \leqslant y_{tn}^\varrho + x_{tn}^\varrho \tag{7r}$$

$$\tfrac{1}{2}\left(\overline{\vartheta}_n - \underline{\vartheta}_n\right) \geqslant x_{tn}^\vartheta \geqslant 0,\ y_{tn}^\vartheta \geqslant 0 \tag{7s}$$

$$\tfrac{1}{2}\left(\overline{\kappa}_\ell - \underline{\kappa}_\ell\right) \geqslant x_{t\ell}^\kappa \geqslant 0,\ y_{t\ell}^\kappa \geqslant 0 \tag{7t}$$

$$\tfrac{1}{2}\left(\overline{\varrho}_n - \underline{\varrho}_n\right) \geqslant x_{tn}^\varrho \geqslant 0,\ y_{tn}^\varrho \geqslant 0 \tag{7u}$$

$$\forall t \in \mathcal{T},\ \forall n \in \mathcal{N},\ \forall \ell \in \mathcal{E},\ ^*\forall \ell \in \mathcal{E}_a \tag{7v}$$

which optimizes the trade-off between the expected cost and pressure variability measure, by choosing only one binary variable (topology) from $v_1,\ldots,v_C$, and by optimizing the multi-stage linear decision rules under the chosen topology.


\end{document}